\documentclass[11pt,a4paper]{article}
\usepackage[margin=30mm]{geometry}
\usepackage{amsmath,amssymb,amsthm,mathtools}
\mathtoolsset{showonlyrefs=true}
\usepackage{microtype}
\usepackage{mathrsfs}
\usepackage[T1]{fontenc}
\usepackage{lmodern}
\usepackage{enumitem}
\usepackage{cite}
\usepackage[
    colorlinks=true,
    linkcolor=blue,
    citecolor=blue,
    urlcolor=blue
]{hyperref}
\numberwithin{equation}{section}
\newtheorem{theorem}{Theorem}[section]
\newtheorem{lemma}[theorem]{Lemma}
\newtheorem{proposition}[theorem]{Proposition}

\theoremstyle{definition}
\newtheorem{definition}[theorem]{Definition}
\theoremstyle{remark}
\newtheorem{remark}[theorem]{Remark}

\newcommand{\dd}{\,\mathrm{d}}
\newcommand{\R}{\mathbb{R}}
\newcommand{\Cc}{C_c^\infty((0,\infty))}

\newcommand{\Del}{\Delta_r}

\DeclareMathOperator{\supp}{supp}

\theoremstyle{plain}
\newtheorem{theoremA}{Theorem}

\begin{document}

\title{Uniqueness of Axisymmetric Vortices \\ in the Magnetic Ginzburg--Landau Equations}

\author{
      {\Large Kazuya Murakami}\thanks{
            Department of Physics, Keio University, Yokohama 223-8533, Japan. 
            \newline\hspace*{1.5em}
            E-mail: \texttt{kazuyamurakami@keio.jp}
      }
}

\thispagestyle{empty}
\maketitle
\begin{abstract}
      We study vortex solutions of the magnetic Ginzburg--Landau model (or the Abelian--Higgs model).
      We establish uniqueness of axisymmetric vortex solutions for every winding number \(n\geq1\) and coupling constant \(\lambda>0\).
      Our proof uses a lattice inequality and constrained minimization to construct a positive radial solution lying above any two given radial solutions.
      A comparison argument based on the virial identity then establishes uniqueness.
      We further prove that the unique radial solution depends real analytically on \(\lambda\) and that both radial profiles are strictly increasing with respect to \(\lambda\).
\end{abstract}

\section{Introduction and main result}

The Abelian--Higgs model, also known as the magnetic Ginzburg--Landau model, arises in the study of superconductivity~\cite{Ginzburg:1950sr,Abrikosov:1957wnz} and Abelian gauge theories~\cite{Nielsen:1973cs,deVega:1976xbp}.
In two spatial dimensions, the energy functional takes the form
\begin{align}\label{eq:AHenergy}
      \mathcal E(\psi,A) &= \frac{1}{2}\int_{\R^2} \left\{|\nabla_A\psi|^2 + |\operatorname{curl}A|^2 + \frac{\lambda}{4}(1 - |\psi|^2)^2\right\} \dd x,
\end{align}
where \(\lambda>0\) is the coupling constant, \(\psi:\R^2\to\mathbb C\) is a complex scalar field, and \(A:\R^2\to\R^2\) is a vector potential.
Here \(\nabla_A := \nabla-iA\) is the covariant derivative, and \(\operatorname{curl}A := \partial_1A_2-\partial_2A_1\) is the magnetic field.
The energy is invariant under gauge transformations
\begin{align}
      (\psi,A) \mapsto (e^{i\chi}\psi, A+\nabla\chi),
\end{align}
where \(\chi\) is a smooth real-valued function.
Critical points of \(\mathcal E\) correspond to static solutions of the field equations and are therefore of physical interest.
The constant configurations \((\psi,A) = (e^{i\alpha},0)\), where \(\alpha\in\R\), are trivial critical points with zero energy.
The model also admits nontrivial finite-energy critical points, called \emph{vortex solutions}.
In this paper we study vortices that are axisymmetric about the origin.

The nonlinear potential term in \eqref{eq:AHenergy} vanishes precisely when \(|\psi|=1\).
Accordingly, we consider finite-energy vortex configurations satisfying \(|\psi(x)|\to1\) as \(|x|\to\infty\).
For such configurations, the scalar field is nonzero outside a sufficiently large disk.
The phase of \(\psi\) therefore defines a continuous map from a large circle \(S_R := \{x\in\R^2:|x|=R\}\) to the unit circle, and its degree
\begin{align}
      n &:= \deg\left(\left. \frac{\psi}{|\psi|} \right|_{S_R} \right) \in \mathbb Z
\end{align}
is independent of sufficiently large \(R\).
This integer is called the \emph{winding number} or \emph{vortex number}.
We shall consider positive winding numbers \(n\geq1\).
The winding number also determines the total magnetic flux, as will be seen explicitly below for axisymmetric vortices.

The case \(\lambda=1\), known as the \emph{critical coupling}, has been extensively studied.
For a smooth finite-energy configuration with winding number \(n\geq1\), the energy admits the Bogomolny decomposition~\cite{Bogomolny:1975de}
\begin{align}
      \mathcal E(\psi,A)
      &= \pi n + \frac{1}{2}\int_{\R^2} \bigg\{
            |(D_1 + i D_2)\psi|^2 + \left(\operatorname{curl}A - \frac{1-|\psi|^2}{2}\right)^2
      \bigg\}\dd x.
\end{align}
Here \(D_j := \partial_j - i A_j\) for \(j=1,2\) and \(i = \sqrt{-1}\).
In particular, \(\mathcal E(\psi,A)\geq\pi n\), and equality holds precisely when
\begin{align}
      (D_1 + i D_2)\psi &= 0,\\
      \operatorname{curl} A &= \frac{1}{2}(1-|\psi|^2). 
\end{align}
These first-order equations are called the \emph{Bogomolny equations}.
Any solution of these equations minimizes the energy within its winding-number class and hence is a critical point of \(\mathcal E\).
Taubes~\cite{Taubes:1979tm} proved that, for any prescribed collection of \(n\) points in \(\R^2\), counted with multiplicity, there exists a unique solution of the Bogomolny equations up to gauge transformations, whose scalar field vanishes precisely at those points.
Moreover, Taubes~\cite{MR581946} showed that every smooth finite-energy solution of the second-order Euler--Lagrange equations with positive winding number at \(\lambda=1\) satisfies the Bogomolny equations; see also Jaffe and Taubes~\cite{Jaffe:1980mj}.

For \(\lambda\neq1\), the above first-order reduction is not available, and the existence and uniqueness questions must be studied directly through the second-order equations.
For general coupling constants \(\lambda>0\), particular attention has been paid to \emph{axisymmetric vortices}, which are of the form
\begin{align}\label{eq:ansatz}
      \psi(r,\theta) = f(r)e^{in\theta},\qquad
      A(r,\theta) = \frac{na(r)}{r}\mathbf e_\theta,
\end{align}
where \(n\geq1\) is the winding number, \((r,\theta)\) are polar coordinates, \(\mathbf e_\theta := (-\sin\theta,\cos\theta)\), and \(f,a\) are real-valued functions of \(r\).
Motivated by regularity at the origin and finite energy, we impose the boundary conditions
\begin{align}\label{eq:bc}
    \lim_{r\downarrow0}(f(r),a(r)) = (0,0),\qquad
    \lim_{r\to\infty}(f(r),a(r)) = (1,1).
\end{align}
For later use, we set
\begin{align}
      P(r) := \frac{1-a(r)}{r},\qquad
      B(r) := (\operatorname{curl}A)(r,\theta)=\frac{na'(r)}{r}. 
\end{align}
Thus \(B\) is the magnetic field under the ansatz \eqref{eq:ansatz}. The boundary conditions also give the magnetic flux of an axisymmetric vortex:
\begin{align}
      \int_{\R^2}\operatorname{curl}A\,\dd x
      = 2\pi\int_0^\infty B(r)r\,\dd r
      = 2\pi n\int_0^\infty a'(r)\,\dd r
      = 2\pi n. 
\end{align}

Substituting \eqref{eq:ansatz} into \eqref{eq:AHenergy} gives \(\mathcal E(\psi,A) = 2\pi E_\lambda(f,a)\), where
\begin{align}\label{eq:radialenergy}
      E_\lambda(f,a) &:= \int_0^\infty\left\{
            \frac{1}{2}f'^2 + \frac{n^2}{2}P^2f^2 + \frac{1}{2}B^2 + \frac{\lambda}{8}(1-f^2)^2
      \right\} r\dd r.
\end{align}
Substituting the same ansatz into the Euler--Lagrange equations of \eqref{eq:AHenergy} yields
\begin{equation}
      \begin{aligned}\label{eq:EL}
            \Del f &= \varphi f,\qquad
            \varphi := n^2 P^2 - \frac{\lambda}{2}(1-f^2),\\
            B' &= -nf^2P.
      \end{aligned}
\end{equation}
Here \(\Del := \frac{\dd^2}{\dd r^2} + \frac{1}{r}\frac{\dd}{\dd r}\) denotes the radial Laplacian.
These are also the Euler--Lagrange equations of the reduced functional \eqref{eq:radialenergy}.

\begin{definition}[Radial solutions]\label{Def:radial}
      A pair \((f,a)\in C^2((0,\infty))^2\) is called a \emph{radial solution} if it satisfies \eqref{eq:EL} and the boundary conditions \eqref{eq:bc}.
      A radial solution is called \emph{positive} if \(f(r)>0\) for all \(r>0\).
\end{definition}

The existence and qualitative behavior of axisymmetric vortices have been studied for general positive coupling constants.
Plohr~\cite{MR641458} obtained detailed asymptotic formulas for vortex profiles at spatial infinity.
Berger and Chen~\cite{MR987294} established the existence of symmetric vortices for arbitrary positive coupling and studied the large-coupling limit.
In particular, a positive radial solution exists for every \(n\geq1\) and \(\lambda>0\).

Uniqueness, however, has been established only in certain ranges of the coupling constant.
Beyond the critical-coupling case discussed above, Clemons~\cite{MR1710018} proved uniqueness for \(\lambda\) sufficiently close to \(1\).
Alama, Bronsard, and Giorgi~\cite{MR1716202} subsequently proved uniqueness of positive axisymmetric solutions under the condition \(\lambda\geq4n^2\).
Their proof uses constrained minimization and a mountain-pass argument.
These results leave the uniqueness question unresolved for general \(\lambda>0\).

The purpose of this paper is to prove uniqueness for the radial Euler--Lagrange equations \eqref{eq:EL} under the boundary conditions \eqref{eq:bc}.
Our main result is as follows.
\begin{theoremA}[Uniqueness]\label{thm:A}
      For every \(n\geq1\) and \(\lambda>0\), \eqref{eq:EL} has a unique radial solution \((f,a)\) satisfying \eqref{eq:bc}.
      Moreover, \(f>0\) in \((0,\infty)\).
\end{theoremA}

\begin{remark}
      The uniqueness assertion of Theorem~\ref{thm:A} also holds for pairs \((f,a)\in H^1_{\mathrm{loc}}((0,\infty))^2\) satisfying \eqref{eq:EL} in the weak sense and \eqref{eq:bc} as limits of their continuous representatives.
      Indeed, standard regularity theory for ODEs implies that
      \((f,a)\in C^\infty((0,\infty))^2\).
      Thus every such weak solution is a radial solution in the sense of Theorem~\ref{thm:A}.
\end{remark}

Theorem~\ref{thm:A} extends the previous uniqueness results to all \(\lambda>0\) without assuming that \(f\) is nonnegative.
Using a lattice inequality (Lemma~\ref{lem:lattice}) and constrained minimization, we construct a positive radial solution lying above any two given radial solutions (Lemma~\ref{lem:upper}).
A comparison argument based on the virial identity (Lemma~\ref{lem:ends}(iii)) then establishes uniqueness (Lemma~\ref{lem:compare}).
We also study the dependence of the unique radial solution on the coupling constant \(\lambda\).
We show that the radial solution depends real analytically on \(\lambda\) and that both components are strictly increasing with respect to \(\lambda\) (Proposition~\ref{prop:analytic}).
The analyticity proof relies on the positivity of the radial Hessian (Proposition~\ref{prop:hessian}) and the analytic implicit-function theorem.
The positivity of the radial Hessian was previously established by Gustafson and Sigal~\cite{MR1772246}, whose proof uses derivatives of the vortex profiles with respect to \(\lambda\).
In contrast, our proof uses a positive vector constructed from the radial derivatives, without assuming differentiability of the radial solution with respect to \(\lambda\).
Thus Proposition~\ref{prop:hessian} provides an alternative proof of their result.

The rest of this paper is organized as follows.
In Section~2, we establish basic properties of radial solutions.
In Section~3, we prove Theorem~\ref{thm:A}.
In Section~4, we establish the positivity of the radial Hessian.
As a consequence of Theorem~\ref{thm:A}, we then prove the real analyticity and strict monotonicity of the radial solution with respect to \(\lambda\).

\section{Properties of radial solutions}
In this section, we prove some properties of radial solutions (see Definition \ref{Def:radial}). 
\begin{lemma}\label{lem:shape}
      Let \((f,a)\) be any radial solution. Then the following properties hold. 
      \begin{enumerate}[label=\textup{(\roman*)},leftmargin=2em]
            \item For any \(r\in(0,\infty)\), we have
            \begin{align}
                  0<a(r)<1,\qquad a'(r)>0,\qquad
                  P(r)>0,\qquad P'(r)<0,\qquad B(r)>0.
            \end{align}
            Furthermore, \(B\) is strictly decreasing and
            \begin{align}
                  \lim_{r\to\infty}B(r) = 0.
            \end{align}
            \item For all \(r \in (0,\infty)\), \(|f(r)|<1\). 
            \item Every zero of \(f\) is simple, and \(f\) has only finitely many zeros.
            \item If \(f\geq0\) in \((0,\infty)\), then \(0<f<1\) and \(f'>0\) in \((0,\infty)\).
      \end{enumerate}
\end{lemma}
\begin{proof}
      We first observe that every zero of \(f\) is simple and isolated.
      Indeed, \(f\) solves the following ODE:
      \begin{align}
            f'' + \frac{1}{r} f' = \left[ n^2 P^2(r) - \frac{\lambda}{2} (1-f^2) \right] f \quad \text{in } (0,\infty).
      \end{align}
      Hence, if \(f(r_0) = f'(r_0) = 0\) for some \(r_0>0\), then the unique solvability of initial-value problems of ODE implies \(f \equiv 0\) in \((0,\infty)\). 
      However, this contradicts \(f(r) \to 1\) as \(r \to \infty\). 
      Thus, every zero of \(f\) is simple and isolated. 

      (i)
      Put \(w:=1-a\).
      Since \(w'/r = -B/n\), the second equation in \eqref{eq:EL} and the boundary conditions \eqref{eq:bc} give
      \begin{align}
            \left(\frac{w'}{r}\right)' = \frac{f^2 w}{r}, \qquad
            w(0^+) = 1, \quad \lim_{r \to \infty} w(r) = 0.
            \label{eq:w}
      \end{align}
      To prove \(w \geq 0\) in \((0,\infty)\), suppose on the contrary that \(w(r_0) < 0\) holds at some \(r_0 \in (0,\infty)\).
      From \(w(0^+) = 1\) and \(w(r) \to 0\) as \(r \to \infty\), we may assume that \(r_0\) is a global minimum point of \(w\) in \((0,\infty)\).
      Since \(f^2(r) w(r) \leq 0 \) as long as \(r>r_0\) and \(w(r) \leq 0\), we infer from \(w'(r_0) = 0\) and \eqref{eq:w} that \(w'\leq 0\) in \((r_0,\infty)\). 
      However, this is not possible due to \(w(r_0) < 0\) and \(w(r) \to 0\) as \(r \to \infty\). 
      Thus, \(w \geq 0\) in \((0,\infty)\).

      Next, we claim \(w>0\) in \((0,\infty)\).
      Again we argue by contradiction and assume that \(w(r_0) = 0 \) holds at some \(r_0 \in (0,\infty)\).
      From \(w \geq 0\), we have \(w'(r_0) = 0\).
      Thus, uniqueness for the initial-value problem associated with \eqref{eq:w} implies \(w \equiv 0\) in \((0,\infty)\). 
      This contradicts \(w(0^+) = 1\). Hence, \(w>0\) in \((0,\infty)\), which is equivalent to \(a<1\) in \((0,\infty)\). 
      
      To prove \(a'>0\), we first observe that \eqref{eq:w} and \(a<1\) in \((0,\infty)\) imply
      \begin{align}
            \left(\frac{a'}{r}\right)' = -\frac{f^2(1-a)}{r} \leq 0 \quad \text{in } (0,\infty). 
      \end{align}
      Since the zeros of \(f\) are isolated and \(1-a > 0\), the function \(a'/r\) is strictly decreasing in \((0,\infty)\). 
      Let \(L := \lim_{r\to\infty} a'(r)/r\in[-\infty, \infty)\).
      If \(L\neq0\), then there exist \(c>0\) and \(R>0\) such that either \(a'(r)\geq cr\) or \(a'(r) \leq - cr\) for all \(r \geq R\).
      In either case, integrating the inequality yields
      \begin{align}
            |a(r) - a(R)| \geq \frac{c}{2}(r^2 - R^2) \quad \text{for all} \ r \geq R, 
      \end{align}
      which contradicts \(a(r)\to1\) as \(r \to \infty\). 
      Thus \(L=0\) holds, and the strict monotonicity of \(a'/r\) implies \(a'/r > 0\) in \((0,\infty)\).
      Consequently, \(a'>0\) in \((0,\infty)\) and  \(B(r) = na'(r)/r \to 0 \) as \(r \to \infty\). 
      In addition, the boundary condition \(a(0^+)=0\) gives \(0<a<1\) in \((0,\infty)\).
      Since \(P=(1-a)/r\), we have \(P>0\) and
      \begin{align}
            P'=-\frac{a'}{r}-\frac{1-a}{r^2}<0
            \qquad\text{in }(0,\infty).
      \end{align}
      Moreover, \(B = na'/r > 0\) in \((0,\infty)\) and \(B\) is strictly decreasing in \((0,\infty)\).
      Thus, the assertions in (i) hold. 

      (ii)
      We argue by contradiction and suppose that \(f(r_0) \geq 1\) holds at some \(r_0 \in (0,\infty)\). 
      From \(f(0^+) = 0\) and \(f(r) \to 1\) as \(r \to \infty\), there exists a global maximum point \(r_1 \in (0,\infty)\) of \(f\) in \((0,\infty)\). 
      Then \eqref{eq:EL} gives
      \begin{align}
            0 \geq \Del f(r_1) = \left(n^2 P^2(r_1) + \frac{\lambda}{2}(f^2(r_1) - 1)\right)f(r_1) > 0,
      \end{align}
      which is a contradiction. Hence \(f < 1\) in \((0,\infty)\). 
      Similarly, suppose that \(f(r_0)\leq-1\) for some \(r_0\in(0,\infty)\).
      Then \(f\) attains a global minimum at some \(r_2\in(0,\infty)\) with \(f(r_2)\leq-1\).
      At this point, \eqref{eq:EL} gives
      \begin{align}
            0 \leq \Del f(r_2) = \left(n^2 P^2(r_2) + \frac{\lambda}{2}(f^2(r_2) - 1)\right)f(r_2) < 0,
      \end{align}
      which is a contradiction. Hence, \(|f| < 1\) in \((0,\infty)\).

      (iii) 
      From the observation at the beginning of the proof, it remains to show that \(f\) has only finitely many zeros. 
      Since \(a(0^+) = 0\), we may choose \(r_0 > 0\) such that \(a(r) \leq 1/2\) for \(0 < r \leq r_0\).
      Set
      \begin{align}
            r_1 := \min\left\{r_0,\frac{n}{2\sqrt{\lambda}}\right\}.
      \end{align}
      Then, using \(|f|<1\) in \((0,\infty)\), we have
      \begin{align}\label{eq:phiorigin}
            \varphi(r)
            \geq \frac{n^2}{4r^2} - \frac{\lambda}{2}
            \geq \frac{n^2}{8r^2} > 0 \quad \text{for every \(r \in (0,r_1]\)}. 
      \end{align}
      Suppose \(f(\rho) = 0\) for some \(\rho\in(0, r_1]\).
      Since \(f(0^+) = f(\rho)=0\) and \(f\not\equiv0\) in \((0,\rho)\), the function \(|f|\) attains a positive maximum at some \(s\in(0,\rho)\).
      Notice that \((-f,a)\) also satisfies \eqref{eq:EL}, since both equations are invariant under the replacement \(f \mapsto -f\).
      Since \(f(s)\neq0\), the function \(|f|\) is smooth near \(s\) and satisfies \(\Del|f|=\varphi|f|\) there.
      At \(s\), we obtain the following contradiction:
      \begin{align}
            0 \geq \Del|f|(s) = \varphi(s)|f(s)| > 0.
      \end{align}
      Thus \(f\) has no zeros in \((0,r_1]\).
      Since \(f(r)\to1\) as \(r \to \infty\) and its zeros are isolated, there are only finitely many zeros.

      (iv)
      Suppose \(f \geq 0\) in \((0,\infty)\). 
      If \(f(r_0) = 0\) holds for some \(r_0 \in (0,\infty)\), then \(f'(r_0) = 0\) and the unique solvability yields \(f \equiv 0\) in \((0,\infty)\). 
      However, this contradicts \(f(r) \to 1\) as \(r \to \infty\). 
      This fact together with (ii) leads to \(0<f<1\) in \((0,\infty)\).

      To show \(f' \geq 0\) in \((0,\infty)\), suppose that \(f'(r_1) < 0\) holds at some \(r_1 \in (0,\infty)\). 
      From \(f(0^+) = 0\), \(f(r) \to 1\) as \(r \to \infty\) and \(0<f<1\) in \((0,\infty)\), 
      we may find $0<r_0 < r_1 <r_2$ such that 
      \(r_0\) (resp. \(r_2\)) is a local maximum (resp. local minimum) of \(f\) with 
      \(f(r_2) < f(r_0)\). From \eqref{eq:EL}, we infer that 
      \[
      \varphi(r_0) f(r_0) = \Del f(r_0) \leq 0 \leq \Del f(r_2) = \varphi(r_2) f(r_2),
      \]
      which gives \(\varphi(r_0) \leq 0 \leq \varphi(r_2)\). 
      On the other hand, the monotonicity of \(a\) in (i) and \(f(r_0) > f(r_2)\) imply 
      \begin{align}
            \varphi(r_2) - \varphi(r_0) = n^2\bigl(P(r_2)^2 - P(r_0)^2\bigr) + \frac{\lambda}{2}\bigl(f(r_2)^2 - f(r_0)^2\bigr) < 0.
      \end{align}
      This is a contradiction and \(f' \geq 0\) in \((0,\infty)\). 

      Finally, suppose that \(f'(r_0) = 0\) for some \(r_0 \in (0,\infty)\).
      Since \(f' \geq 0\), we have \(f''(r_0) = 0\). 
      By \(f(r_0) > 0\) and \(0= \Del f(r_0) = \varphi(r_0) f(r_0) \), we observe \(\varphi(r_0) = 0\).
      Moreover, (i) gives 
      \begin{align}
            \varphi'(r_0) = 2 n^2 P(r_0) P'(r_0) < 0.
      \end{align}
      Hence there exists \(\varepsilon > 0\) such that \(\varphi < 0\) in \((r_0,r_0+\varepsilon)\).
      Since \((rf')' = r\varphi f\), we obtain
      \begin{align}
      rf'(r)
      =\int_{r_0}^r s\varphi(s)f(s)\,ds<0
      \quad \text{for every \(r \in (r_0,r_0+\varepsilon)\)},
      \end{align}
      which contradicts \(f' \geq 0\) in \((0,\infty)\). 
      Therefore \(f' > 0\) in \((0, \infty)\).
\end{proof}

We next establish the endpoint estimates needed in the subsequent analysis.
For more precise asymptotic formulas at infinity, see Plohr~\cite{MR641458}.
\begin{lemma}\label{lem:ends}
      Let \((f,a)\) be any radial solution. Then the following properties hold.
      \begin{enumerate}[label=\textup{(\roman*)},leftmargin=2em]
            \item The limit \(B(0) := \lim_{r\downarrow0}B(r)\) exists in \(\mathbb R\), and \(rf'(r)\to0\) as \(r\downarrow0\).
            Moreover,
            \begin{align}
                  a(r) = O(r^2), \qquad
                  |f(r)| = O(r^n), \qquad
                  |f'(r)| = O(r^{n-1})
                  \qquad (r\downarrow0).
            \end{align}
            \item There exists \(c>0\) such that
            \begin{align}
                  |1-f(r)| + |f'(r)| + |1-a(r)| + |B(r)| = O(e^{-cr})
                  \qquad (r\to\infty).
            \end{align}
            \item The following identity holds:
            \begin{align}\label{eq:virial}
                  \int_0^\infty B^2r\dd r
                  = \frac{\lambda}{4}\int_0^\infty (1-f^2)^2r\dd r.
            \end{align}
      \end{enumerate}
\end{lemma}
\begin{proof}
      (i)
      As shown in the proof of Lemma~\ref{lem:shape}(iii), the function \(f\) has no zeros in \((0,r_1]\), where \(r_1>0\) is chosen in \eqref{eq:phiorigin}.
      Hence, \(f\) has constant sign in \((0,r_1]\).
      Since \((-f,a)\) also satisfies \eqref{eq:EL}, the function \(|f|\) satisfies \(\Del|f|=\varphi|f|\) in \((0,r_1]\).
      By \eqref{eq:EL} and \eqref{eq:phiorigin}, we have
      \begin{align}
            (r|f|'(r))' = r\varphi(r)|f(r)| > 0
            \qquad\text{for all }r \in (0,r_1].
      \end{align}
      Thus, \(r|f|'(r)\) is strictly increasing in \((0,r_1]\), and the limit \(L:=\lim_{r\downarrow0}r|f|'(r) \in [-\infty,\infty)\) exists.
      Suppose that \(L\neq0\).
      Then there exist \(c>0\) and \(R\in(0,r_1]\) such that either \(r|f|'(r) \geq c\) for all \(r\in(0,R)\), or \(r|f|'(r) \leq -c\) for all \(r\in(0,R)\).
      In either case, integration gives
      \begin{align}
            \bigl||f(R)| - |f(r)|\bigr| &= \left|\int_r^R |f|'(s)\,\dd s\right|\\
            &\geq c\int_r^R \frac{\dd s}{s} = c\log\frac{R}{r}
            \qquad\text{for all }r\in(0,R).
      \end{align}
      Letting \(r\downarrow0\), we obtain a contradiction to \(|f|<1\) in Lemma~\ref{lem:shape}(ii).
      Therefore \(L=0\).
      Since \(f\) has constant sign in \((0,r_1]\), we conclude that \(rf'(r)\to0\) as \(r\downarrow0\).

      Next, set \(\gamma:=n/(2\sqrt2)>0\).
      By \eqref{eq:EL} and \eqref{eq:phiorigin}, we have
      \begin{align}
            \Del|f| \geq \frac{\gamma^2}{r^2}|f|
            \qquad\text{in }(0,r_1].
      \end{align}
      On the other hand, a direct computation gives
      \begin{align}
            \Del r^\gamma = \frac{\gamma^2}{r^2}r^\gamma.
      \end{align}
      Put
      \begin{align}
            h(r) := |f(r)| - |f(r_1)|\left(\frac{r}{r_1}\right)^\gamma.
      \end{align}
      Combining the two preceding relations, we obtain
      \begin{align}
            \Del h\geq\frac{\gamma^2}{r^2}h
            \qquad\text{in }(0,r_1].
      \end{align}
      Moreover, \(h(r_1)=0\), and \(h(r)\to0\) as \(r\downarrow0\) by the boundary condition \eqref{eq:bc}.
      Suppose that \(h(r)>0\) for some \(r\in(0,r_1)\).
      Then \(h\) attains a positive maximum at some \(s\in(0,r_1)\).
      At \(s\), we obtain the following contradiction:
      \begin{align}
            0 \geq \Del h(s) \geq \frac{\gamma^2}{s^2}h(s) > 0.
      \end{align}
      Thus \(h\leq0\) in \((0,r_1]\), which gives
      \begin{align}
            |f(r)| \leq |f(r_1)|\left(\frac{r}{r_1}\right)^\gamma
            \qquad\text{for all }r \in (0,r_1].
      \end{align}
      In particular, \(f(r) = O(r^\gamma)\) as \(r\downarrow0\).

      Since \(0<a<1\) by Lemma~\ref{lem:shape}(i), the second equation in \eqref{eq:EL} gives
      \begin{align}
            |B'(r)| = \frac{nf(r)^2(1-a(r))}{r}
            \leq \frac{nf(r)^2}{r} = O(r^{2\gamma-1})
            \qquad(r\downarrow0).
      \end{align}
      Since \(\gamma>0\), the function \(B'\) is integrable near the origin.
      Hence, the limit \(B(0)\) exists and is finite.
      Using \(a'(r)=rB(r)/n\) and \(a(0^+)=0\) from \eqref{eq:bc}, we obtain
      \begin{align}
            a(r) = \frac{1}{n}\int_0^r sB(s)\,\dd s = O(r^2)
            \qquad(r\downarrow0).
      \end{align}
      In particular, there exist \(C_1>0\) and \(\rho_0 > 0\) such that
      \begin{align}
            \varphi(r) &= \frac{n^2(1-a(r))^2}{r^2} - \frac{\lambda}{2}(1-f(r)^2)\\
            &\geq \frac{n^2}{r^2} - C_1
            \qquad\text{for all }r \in (0,\rho_0].
      \end{align}
      Here we used \((1-a)^2 \geq 1-2a\), \(a(r) = O(r^2)\), and \(|f|<1\).

      Set \(\beta := C_1/(2(n+1))\), and put \(\zeta(r) := r^n e^{-\beta r^2}\).
      A direct calculation gives
      \begin{align}
            \Del \zeta(r) = \left(\frac{n^2}{r^2} - 2C_1 + 4\beta^2r^2\right)\zeta(r).
      \end{align}
      Choose \(\rho\in(0,\min\{\rho_0,r_1\}]\) sufficiently small
      that
      \begin{align}
            4\beta^2 r^2 \leq C_1, \qquad
            \frac{n^2}{r^2} > C_1
            \qquad\text{for all }r\in(0,\rho].
      \end{align}
      Then
      \begin{align}
            \Del \zeta(r) \leq \left(\frac{n^2}{r^2} - C_1\right)\zeta(r)
            \qquad\text{for all }r\in(0,\rho].
      \end{align}
      Put
      \begin{align}
            k(r) := |f(r)| - \frac{\zeta(r)}{\zeta(\rho)}.
      \end{align}
      Since \(\Del|f| = \varphi|f|\), the preceding inequalities give
      \begin{align}
            \Del k \geq \left(\frac{n^2}{r^2} - C_1\right)k
            \qquad\text{in }(0,\rho].
      \end{align}
      Moreover, \(k(\rho) = |f(\rho)|-1 < 0\) by Lemma~\ref{lem:shape}(ii), and \(k(r)\to0\) as \(r\downarrow0\).
      Suppose that \(k(r)>0\) for some \(r\in(0,\rho)\).
      Then \(k\) attains a positive maximum at some \(s\in(0,\rho)\).
      At \(s\), we obtain
      \begin{align}
            0 \geq \Del k(s) \geq \left(\frac{n^2}{s^2} - C_1\right)k(s) > 0,
      \end{align}
      which is a contradiction.
      Thus \(k \leq 0\) in \((0,\rho]\), and hence
      \begin{align}
            |f(r)| \leq \frac{r^ne^{-\beta r^2}}{\rho^ne^{-\beta\rho^2}}
            \qquad\text{for all }r\in(0,\rho].
      \end{align}
      This proves \(|f(r)| = O(r^n)\) as \(r\downarrow0\).

      Finally, the first equation in \eqref{eq:EL} gives
      \begin{align}
            (rf'(r))' = r\varphi(r)f(r).
      \end{align}
      Integrating this identity over \((\varepsilon,r)\) and using \(\lim_{\varepsilon\downarrow0}\varepsilon f'(\varepsilon) = 0\), we obtain
      \begin{align}
            rf'(r) = \int_0^r s\varphi(s)f(s)\,\dd s.
      \end{align}
      Since \(0<a<1\) and \(|f|<1\), we have
      \begin{align}
            |\varphi(r)| \leq \frac{n^2}{r^2} + \frac{\lambda}{2}
            \qquad\text{for all }r \in (0,\infty).
      \end{align}
      Together with \(|f(r)| = O(r^n)\), this yields
      \begin{align}
            |f'(r)| &\leq \frac{1}{r}\int_0^r s|\varphi(s)f(s)|\,\dd s\\
            &= O(r^{n-1}) \qquad(r\downarrow0).
      \end{align}
      Thus, the assertions in (i) hold.

      (ii)
      By \eqref{eq:bc} and Lemma~\ref{lem:shape}(ii), we may choose \(R \geq 2\) such that
      \begin{align}
            \frac{1}{2} \leq f(r) < 1
            \qquad\text{for all }r \in [R,\infty).
      \end{align}
      Put \(w:=1-a\).
      By Lemma~\ref{lem:shape}(i), we have \(w>0\) in \((0,\infty)\).
      The second equation in \eqref{eq:EL} gives
      \begin{align}
            w'' - \frac{w'}{r} = f^2w \geq \frac{1}{4}w
            \qquad\text{in }[R,\infty).
      \end{align}
      Set \(\xi(r) := w(R)e^{-(r-R)/4}\).
      By the definition of \(\xi\), together with \(R\geq2\) and \(f^2 \geq 1/4\) in \([R,\infty)\), we obtain
      \begin{align}
            \xi'' - \frac{\xi'}{r} - f^2\xi
            &= \left(\frac{1}{16} + \frac{1}{4r} - f^2\right)\xi\\
            &\leq -\frac{1}{16}\xi < 0
            \qquad\text{in }[R,\infty).
      \end{align}
      Thus, the function \(g := w-\xi\) satisfies
      \begin{align}
            g'' - \frac{g'}{r} - f^2g > 0
            \qquad\text{in }(R,\infty).
      \end{align}
      Moreover, \(g(R)=0\), and \(g(r)\to0\) as \(r\to\infty\) by \eqref{eq:bc}.
      Suppose that \(g(r)>0\) for some \(r>R\).
      Then \(g\) attains a positive maximum at some \(s>R\).
      At \(s\), we have
      \begin{align}
            0 < g''(s) - \frac{g'(s)}{s} - f(s)^2 g(s)
            \leq -f(s)^2 g(s) < 0,
      \end{align}
      which is a contradiction.
      Hence,
      \begin{align}
            0 < w(r) \leq w(R)e^{-(r-R)/4}
            \qquad\text{in }[R,\infty).
      \end{align}

      By Lemma~\ref{lem:shape}(i), \(B(r)\to0\) as \(r\to\infty\).
      Integrating the second equation in \eqref{eq:EL}, we obtain
      \begin{align}
            B(r) = n\int_r^\infty \frac{f(s)^2w(s)}{s}\,\dd s.
      \end{align}
      Since \(|f|<1\) by Lemma~\ref{lem:shape}(ii), the above estimate for \(w\) gives
      \begin{align}
            0 < B(r) &\leq n\int_r^\infty \frac{w(s)}{s}\,\dd s\\
            &\leq \frac{nw(R)}{r}\int_r^\infty e^{-(s-R)/4}\,\dd s\\
            &= \frac{4nw(R)}{r} e^{-(r-R)/4}
            \qquad\text{for all }r \in [R,\infty).
      \end{align}
      Consequently, \(1-a\) and \(B\) decay exponentially.
      Since \(P=w/r\), we also have
      \begin{align}
            P(r)^2 = O(e^{-r/2})\qquad(r\to\infty).
      \end{align}

      Next, put \(u := 1-f\).
      By Lemma~\ref{lem:shape}(ii) and the choice of \(R\), we have \(0 < u \leq 1/2\) in \([R,\infty)\).
      Combining the first equation in \eqref{eq:EL} with the above estimate for \(w\), we obtain
      \begin{align}
            \Del u &= \frac{\lambda}{2} f(1+f)u - n^2 P^2 f\\
            &\geq \frac{3\lambda}{8}u - Ce^{-r/2}
            \qquad\text{for all }r \in [R,\infty),
      \end{align}
      for some constant \(C>0\).
      
      Set \(\mu := 3\lambda/8\), and choose \(\kappa \in (0, \min\{1/2, \sqrt\mu\})\).
      Since \(\kappa < 1/2\) and \(\mu - \kappa^2 > 0\), we may choose \(A>0\) sufficiently large that
      \begin{align}
            (\mu - \kappa^2)Ae^{-\kappa r} &\geq Ce^{-r/2} \qquad\text{for all }r \in [R,\infty),\\
            Ae^{-\kappa R} &\geq u(R).
      \end{align}
      On the other hand, a direct computation gives
      \begin{align}
            \Del e^{-\kappa r} = \left(\kappa^2 - \frac{\kappa}{r}\right)e^{-\kappa r}.
      \end{align}
      Put \(q(r) := u(r)-Ae^{-\kappa r}\).
      Then
      \begin{align}
            \Del q &\geq \mu u - Ce^{-r/2} - A\left(\kappa^2 - \frac{\kappa}{r}\right)e^{-\kappa r}\\
            &\geq \mu q
            \qquad\text{for all }r \in [R,\infty).
      \end{align}
      Moreover, \(q(R)\leq0\) by the choice of \(A\), and \(q(r)\to0\) as \(r\to\infty\) by \eqref{eq:bc}.
      Suppose that \(q(r)>0\) for some \(r>R\).
      Then \(q\) attains a positive maximum at some \(s>R\).
      At \(s\), we obtain
      \begin{align}
            0 \geq \Del q(s) \geq \mu q(s)>0,
      \end{align}
      which is a contradiction.
      Therefore,
      \begin{align}
            0 < 1-f(r) \leq Ae^{-\kappa r}
            \qquad\text{for all }r \in [R,\infty).
      \end{align}

      By the first equation in \eqref{eq:EL} and Lemma~\ref{lem:shape}(ii), we have
      \begin{align}
            |\Del u| &= \left|\frac{\lambda}{2}f(1+f)u - n^2P^2f\right|\\
            &\leq \lambda|u| + n^2P^2
            \qquad\text{in }[R,\infty).
      \end{align}
      Since \(u = O(e^{-\kappa r})\), \(P^2 = O(e^{-r/2})\), and \(\kappa<1/2\), we obtain
      \begin{align}
            |\Del u| = O(e^{-\kappa r}) \qquad(r\to\infty).
      \end{align}
      Consequently, the identity
      \begin{align}
            (ru'(r))' = r\Del u(r)
      \end{align}
      shows that \(ru'(r)\) has a finite limit \(L_\infty:=\lim_{r\to\infty}ru'(r)\), because \(r\Del u(r)\) is integrable on \((R,\infty)\).
      Suppose that \(L_\infty\neq0\).
      Then there exist \(c>0\) and \(R_1\geq R\) such that either \(ru'(r)\geq c\) for all \(r\geq R_1\), or \(ru'(r)\leq-c\) for all \(r\geq R_1\).
      In either case,
      \begin{align}
            |u(r) - u(R_1)| &= \left|\int_{R_1}^r u'(s)\,\dd s\right|\\
            &\geq c\log\frac{r}{R_1}
            \qquad\text{for all }r \in [R_1,\infty),
      \end{align}
      which contradicts \(u(r)\to0\) as \(r\to\infty\).
      Thus \(L_\infty=0\), and integration from \(r\) to infinity yields
      \begin{align}
            ru'(r) = -\int_r^\infty s\Del u(s)\,\dd s.
      \end{align}
      Hence
      \begin{align}
            |f'(r)| = |u'(r)|
            &\leq \frac{C}{r}\int_r^\infty se^{-\kappa s}\,\dd s\\
            &= O(e^{-\kappa r})\qquad(r\to\infty).
      \end{align}
      The estimates above prove (ii).

      (iii)
      Define
      \begin{align}
            \mathcal T(r)
            := \frac{1}{2}(rf'(r))^2 - \frac{n^2}{2}(1-a(r))^2f(r)^2 - \frac{\lambda}{8}r^2(1-f(r)^2)^2 + \frac{1}{2}r^2B(r)^2.
      \end{align}
      Differentiating \(\mathcal T\) and using \(a'(r) = rB(r)/n\), we obtain
      \begin{align}
            \mathcal T'(r) &= r^2f'(r)\bigl(\Del f(r)-\varphi(r)f(r)\bigr) + r^2B(r)\bigl(B'(r)+nf(r)^2P(r)\bigr)\\
            &\quad + r\left\{B(r)^2 - \frac{\lambda}{4}(1-f(r)^2)^2\right\}.
      \end{align}
      By \eqref{eq:EL}, the first two terms on the right-hand side vanish.
      Therefore,
      \begin{align}\label{eq:virialT}
            \mathcal T'(r) = r\left\{B(r)^2 - \frac{\lambda}{4}(1-f(r)^2)^2\right\}
            \qquad\text{for all }r \in (0,\infty).
      \end{align}
      
      By (i), we have \(rf'(r)\to0\), \(f(r)\to0\), \(a(r)\to0\), and \(B(r) = O(1)\) as \(r\downarrow0\).
      On the other hand, by (ii), the functions \(f'\), \(1-a\), \(1-f\), and \(B\) decay exponentially
      as \(r\to\infty\).
      These estimates imply
      \begin{align}
            \lim_{r\downarrow0}\mathcal T(r)
            = \lim_{r\to\infty}\mathcal T(r)
            = 0.
      \end{align}
      Integrating \eqref{eq:virialT} over \((0,\infty)\) and using the above limits, we obtain
      \begin{align}
            \int_0^\infty B(r)^2r\,\dd r
            = \frac{\lambda}{4}\int_0^\infty(1-f(r)^2)^2r\,\dd r,
      \end{align}
      which proves \eqref{eq:virial} and completes the proof.
\end{proof}

\section{Proof of uniqueness}
For pairs \(u=(g,b)\) and \(\widetilde u=(\widetilde g,\widetilde b)\), let \(u\vee\widetilde u\) and \(u\wedge\widetilde u\) denote their componentwise maximum and minimum, respectively.

\begin{lemma}\label{lem:lattice}
      The radial energy \(E_\lambda\) has the following properties.
      \begin{enumerate}[label=\textup{(\roman*)},leftmargin=2em]
            \item For any pair \((g,b)\in AC_{\mathrm{loc}}(0,\infty)^2\) of real-valued profiles, the following identities and inequalities hold:
            \begin{align}
                  E_\lambda(|g|,b) &= E_\lambda(g,b),\\
                  E_\lambda(g,\min\{b,1\}) &\leq E_\lambda(g,b),\\
                  E_\lambda(\min\{g,1\},b) &\leq E_\lambda(g,b)
                  \qquad\text{if }g\geq0\text{ in }(0,\infty).
            \end{align}
            \item For any \(u=(g,b)\) and \(\widetilde u=(\widetilde g,\widetilde b)\) in \(AC_{\mathrm{loc}}(0,\infty)^2\) satisfying \(g, \widetilde g \geq 0\) and \(b, \widetilde b \leq 1\) in \((0,\infty)\), the following inequality holds:
            \begin{align}\label{eq:lattice}
                  E_\lambda(u\vee\widetilde u) + E_\lambda(u\wedge\widetilde u)
                  \leq E_\lambda(u) + E_\lambda(\widetilde u).
            \end{align}
      \end{enumerate}
\end{lemma}
\begin{proof}
      (i)
      By \eqref{eq:radialenergy}, the radial energy of \((g,b)\) can be written as
      \begin{align}
            E_\lambda(g,b)
            = \int_0^\infty \left\{
                  \frac{1}{2}|g'|^2 + \frac{n^2}{2r^2}g^2(1-b)^2 + \frac{n^2}{2r^2}|b'|^2 + \frac{\lambda}{8}(1-g^2)^2
            \right\}r\dd r.
      \end{align}      
      We first consider the replacement of \(g\) by \(|g|\).
      Since \(g\in AC_{\mathrm{loc}}(0,\infty)\), we have
      \begin{align}
            ||g|'| = |g'| \qquad\text{almost everywhere in }(0,\infty).
      \end{align}
      Since the remaining terms in the integrand depend on \(g\) only through \(g^2\), we obtain
      \begin{align}
            E_\lambda(|g|,b) = E_\lambda(g,b).
      \end{align}
      Next, put \(\widehat b := \min\{b,1\}\).
      Since
      \begin{align}
            |\widehat b'| \leq |b'|, \qquad
            |1-\widehat b| \leq |1-b|
            \qquad\text{almost everywhere in }(0,\infty),
      \end{align}
      we obtain \(E_\lambda(g,\widehat b) \leq E_\lambda(g,b)\).

      Finally, suppose that \(g\geq0\) in \((0,\infty)\), and put \(\widehat g := \min\{g,1\}\).
      Then
      \begin{align}
            |\widehat g'|^2 \leq |g'|^2, \qquad
            0 \leq \widehat g^2 \leq g^2, \qquad
            (1-\widehat g^2)^2 \leq (1-g^2)^2
            \quad\text{almost everywhere in }(0,\infty).
      \end{align}
      Since \((1-b)^2 \geq 0\), none of the terms in the integrand increases under this replacement.
      Hence, \(E_\lambda(\widehat g, b) \leq E_\lambda(g,b)\), proving (i).

      (ii)
      We next consider the componentwise maximum and minimum of \(u=(g,b)\) and \(\widetilde u=(\widetilde g,\widetilde b)\).
      For any two locally absolutely continuous functions \(h\) and \(\widetilde h\), the functions \(h\vee\widetilde h\) and \(h\wedge\widetilde h\) are also locally absolutely continuous, and
      \begin{align}
            |(h\vee\widetilde h)'|^2 + |(h\wedge\widetilde h)'|^2
            = |h'|^2 + |\widetilde h'|^2
            \qquad\text{almost everywhere in }(0,\infty).
      \end{align}
      Indeed, where \(h>\widetilde h\) or \(h<\widetilde h\), the maximum and minimum simply interchange the two functions.
      On the set where \(h=\widetilde h\), their derivatives agree almost everywhere, so the identity also holds there.
      Applying this identity to \(g,\widetilde g\) and to \(b,\widetilde b\), we see that the sum of the terms involving \(|g'|^2\) and \(|b'|^2\) is unchanged.
      Similarly, we have
      \begin{align}
            (1 - (g\vee\widetilde g)^2)^2 + (1 - (g\wedge\widetilde g)^2)^2
            = (1-g^2)^2 + (1-\widetilde g^2)^2.
      \end{align}

      It remains to consider the term involving \(g^2(1-b)^2\).
      Write \((g_+,b_+) := u\vee\widetilde u\) and \((g_-,b_-) :=u \wedge\widetilde u\).
      Put
      \begin{align}
            F(g,b) := g^2(1-b)^2.
      \end{align}
      Fix \(r\in(0,\infty)\).
      If \(g(r)\geq\widetilde g(r)\) and \(b(r)\geq\widetilde b(r)\), then \(u\vee\widetilde u=u\) and \(u\wedge\widetilde u=\widetilde u\) at \(r\).
      If \(g(r)\leq\widetilde g(r)\) and \(b(r)\leq\widetilde b(r)\), then \(u\vee\widetilde u=\widetilde u\) and \(u\wedge\widetilde u=u\) at \(r\).
      In either case, the sum of the \(F\)-terms is unchanged.

      We next consider the case where \(g(r)\geq\widetilde g(r)\) and \(b(r)<\widetilde b(r)\).
      Then
      \begin{align}
            F(g,b) + F(\widetilde{g},\widetilde{b}) - F(g,\widetilde b) - F(\widetilde g,b)
            = (g^2 - \widetilde g^2)\bigl((1-b)^2 - (1-\widetilde b)^2\bigr) \geq 0,
      \end{align}
      where all functions are evaluated at \(r\).
      Indeed, \(g \geq \widetilde g \geq 0\) implies \(g^2 - \widetilde g^2 \geq 0\), while \(b < \widetilde b \leq 1\) implies \((1-b)^2 - (1-\widetilde b)^2 \geq 0\).
      The case \(g(r)<\widetilde g(r)\) and \(b(r)\geq\widetilde b(r)\) follows by interchanging the two pairs.
      Thus, in every case,
      \begin{align}
            F(g_+,b_+) + F(g_-,b_-) \leq F(g,b) + F(\widetilde g, \widetilde b).
      \end{align}
      Integrating this inequality, together with the identities for the other terms in the energy, yields \eqref{eq:lattice}.
\end{proof}

\begin{lemma}\label{lem:sub}
      Let \((f,a)\) be any radial solution and put \(v:=(|f|,a)\).
      For any \(\vartheta := (\vartheta_1,\vartheta_2) \in (H_c^1)^2\) satisfying \(\vartheta_1\geq0\) in \((0,\infty)\), the following identity holds:
      \begin{align}\label{eq:sub}
            \dd E_\lambda(v)[\vartheta] = -2\sum_{r_k\in\supp\vartheta_1} r_k|f'(r_k)|\vartheta_1(r_k) \leq 0,
      \end{align}
      where the sum runs over the zeros \(r_k\) of \(f\) in \(\supp\vartheta_1\).
      Here \(H_c^1\) denotes the space of \(H^1((0,\infty))\) functions with compact support in \((0,\infty)\).
      The point values \(\vartheta_1(r_k)\) in \eqref{eq:sub} are understood as the values of the unique continuous representative of \(\vartheta_1\).
      In particular, if \(f>0\) in \((0,\infty)\), then \(\dd E_\lambda(v)[\vartheta] = 0\).
\end{lemma}
\begin{proof}
      By Lemma~\ref{lem:ends}, we have \(E_\lambda(f,a)<\infty\).
      Lemma~\ref{lem:lattice}(i) then gives
      \begin{align}
            E_\lambda(v) = E_\lambda(|f|,a) = E_\lambda(f,a) < \infty.
      \end{align}
      Since \(\vartheta\) has compact support in \((0,\infty)\), we may choose \(0<\alpha<\beta<\infty\) such that \(\supp\vartheta_1\cup\supp\vartheta_2 \subset (\alpha,\beta)\).
      The integrand of \(E_\lambda(v+t\vartheta) - E_\lambda(v)\) vanishes outside \((\alpha,\beta)\).
      Since the energy density in \eqref{eq:radialenergy} is polynomial in the profiles and their first derivatives, the difference \(E_\lambda(v+t\vartheta) - E_\lambda(v)\) is a polynomial in \(t\).
      Thus, the first variation \(\dd E_\lambda(v)[\vartheta]\) is well defined.

      We first consider the variation with respect to \(a\).
      Since \(B = na'/r\) and \(P=(1-a)/r\), differentiating \eqref{eq:radialenergy} gives
      \begin{align}
            \dd E_\lambda(v)[(0,\vartheta_2)]
            &= \int_0^\infty \left\{nB\vartheta_2' - n^2P|f|^2\vartheta_2\right\}\dd r\\
            &= -n\int_0^\infty \left\{B' + nf^2P\right\}\vartheta_2\dd r\\
            &= 0,
      \end{align}
      where we used integration by parts, the compact support of \(\vartheta_2\), and the second equation in \eqref{eq:EL}.

      Next, the variation with respect to the scalar profile is given by
      \begin{align}\label{eq:subscalar}
            \dd E_\lambda(v)[(\vartheta_1,0)]
            = \int_0^\infty \left\{|f|'\vartheta_1' + \varphi|f|\vartheta_1\right\}r\dd r.
      \end{align}
      By Lemma~\ref{lem:shape}(iii), the zeros of \(f\) are simple and finite.
      Let \(r_1<\cdots<r_m\) be the zeros of \(f\) in \((\alpha,\beta)\), and put \(r_0:=\alpha\) and \(r_{m+1}:=\beta\).
      On each interval \((r_j,r_{j+1})\), the function \(f\) has constant sign.
      Since \((-f,a)\) also satisfies \eqref{eq:EL}, the function \(|f|\) is smooth on this interval and satisfies
      \begin{align}
            (r|f|'(r))' = r\varphi(r)|f(r)|
            \qquad\text{for all }r\in(r_j,r_{j+1}).
      \end{align}
      Applying integration by parts to \eqref{eq:subscalar} on each of these intervals, we obtain
      \begin{align}
            \dd E_\lambda(v)[(\vartheta_1,0)]
            &= \sum_{j=0}^{m}\left[\,r|f|'(r)\vartheta_1(r)\,\right]_{r_j+0}^{r_{j+1}-0}\\
            &= \sum_{k=1}^{m}r_k\left\{|f|'(r_k-0)-|f|'(r_k+0)\right\}\vartheta_1(r_k),
      \end{align}
      where the boundary terms at \(\alpha\) and \(\beta\) vanish because \(\vartheta_1\) has compact support in \((\alpha,\beta)\).

      Since every zero \(r_k\) is simple, the sign of \(f\) changes at \(r_k\).
      Consequently, the one-sided derivatives of \(|f|\) satisfy
      \begin{align}
            |f|'(r_k-0) = -|f'(r_k)|,\qquad
            |f|'(r_k+0) = |f'(r_k)|.
      \end{align}
      Therefore,
      \begin{align}
            \dd E_\lambda(v)[\vartheta]
            &= -2\sum_{k=1}^{m}r_k|f'(r_k)|\vartheta_1(r_k)\\
            &= -2\sum_{r_k\in\supp\vartheta_1}r_k|f'(r_k)|\vartheta_1(r_k) \leq 0,
      \end{align}
      where the last inequality follows from \(\vartheta_1\geq0\).
      This proves \eqref{eq:sub}.
      If \(f>0\) in \((0,\infty)\), then \(f\) has no zeros, and hence the first variation vanishes.
\end{proof}

\begin{lemma}\label{lem:upper}
      Let \(u_1,\ldots,u_m\) be radial solutions with \(m\geq1\), and put
      \begin{align}
            u_i := (f_i,a_i), \qquad
            v_i := (|f_i|,a_i), \qquad
            v := v_1\vee\cdots\vee v_m.
      \end{align}
      Then there exists a positive radial solution \(u^\circ=(g^\circ,b^\circ)\) such that
      \begin{align}
            u^\circ \geq v\quad\text{in }(0,\infty).
      \end{align}
      Here, the inequality is understood componentwise.
\end{lemma}
\begin{proof}
      We divide the proof into three steps.
      We first construct a minimizer under the constraint \(v\leq(g,b)\leq\mathbf1\), where \(\mathbf1:=(1,1)\).
      We then prove that this minimizer satisfies the Euler--Lagrange equations.
      Finally, we verify the boundary conditions and positivity.

      \textbf{Step 1.}
      In this step, we prove the existence of a minimizer of \(E_\lambda\) in the class \(\mathcal K\) defined below.

      Define
      \begin{align}
            \mathcal K := \bigl\{
                  (g,b)\in AC_{\mathrm{loc}}(0,\infty)^2:
                  v \leq (g,b) \leq \mathbf1, \quad
                  b(0^+) = 0, \quad E_\lambda(g,b) < \infty
            \bigr\}.
      \end{align}
      We first show that \(\mathcal K\) is nonempty.
      By Lemma~\ref{lem:ends}, each radial solution \(u_i\) has finite energy.
      Lemma~\ref{lem:lattice}(i) therefore gives
      \begin{align}
            E_\lambda(v_i) = E_\lambda(u_i) < \infty \qquad(1\leq i\leq m).
      \end{align}
      By Lemma~\ref{lem:shape}(i)--(ii), we have
      \begin{align}
            0 \leq v_i \leq \mathbf1 \qquad\text{in }(0,\infty).
      \end{align}
      Repeated application of Lemma~\ref{lem:lattice}(ii), together with the nonnegativity of the energy, yields
      \begin{align}
            E_\lambda(v) \leq \sum_{i=1}^m E_\lambda(v_i) < \infty.
      \end{align}
      Moreover, the second component of \(v\) tends to zero as \(r\downarrow0\), since \(a_i(0^+) = 0\) for every \(i\).
      Thus, \(v\in\mathcal K\).

      Let \((g_j,b_j)\subset\mathcal K\) be a minimizing sequence.
      Since \(\inf_{\mathcal K}E_\lambda \leq E_\lambda(v) < \infty\), we may assume that
      \begin{align}
            E_\lambda(g_j,b_j) \leq C\qquad\text{for all }j,
      \end{align}
      where \(C>0\) is independent of \(j\).
      Since all terms in the integrand of \eqref{eq:radialenergy} are nonnegative, the energy bound above implies
      \begin{align}
            \int_0^\infty \frac{|b_j'(r)|^2}{r}\,\dd r
            \leq \frac{2}{n^2}E_\lambda(g_j,b_j)
            \leq \frac{2C}{n^2}.
      \end{align}
      Since \(b_j(0^+)=0\), the Cauchy--Schwarz inequality implies
      \begin{align}
            0\leq b_j(r)
            &= \int_0^r b_j'(s)\,\dd s\\
            &\leq \left(\int_0^r\frac{|b_j'(s)|^2}{s}\,\dd s\right)^{1/2} \left(\int_0^r s\,\dd s\right)^{1/2}\\
            &\leq \frac{\sqrt C}{n}\,r
            =: C_0r \qquad\text{for all }r>0.
      \end{align}

      Fix \(0<\delta<R<\infty\).
      Since all terms in the integrand of \eqref{eq:radialenergy} are nonnegative, the energy bound \(E_\lambda(g_j,b_j)\leq C\) gives
      \begin{align}
            \int_\delta^R |g_j'(r)|^2\,\dd r &\leq \frac{2C}{\delta},\\
            \int_\delta^R |b_j'(r)|^2\,\dd r &\leq \frac{2CR}{n^2}.
      \end{align}
      Moreover, since \(0 \leq g_j, b_j \leq 1\), we have
      \begin{align}
            \int_\delta^R \bigl(|g_j(r)|^2 + |b_j(r)|^2\bigr)\,\dd r
            \leq 2(R-\delta).
      \end{align}
      Combining these estimates, we obtain
      \begin{align}
            \sup_j\bigl(\|g_j\|_{H^1(\delta,R)} + \|b_j\|_{H^1(\delta,R)}\bigr) < \infty.
      \end{align}
      By the Rellich theorem, weak compactness in \(H^1(\delta,R)\), and a diagonal argument, there exist a subsequence and a pair \(u^\circ := (g^\circ,b^\circ)\) such that
      \begin{align}
            (g_j,b_j) &\to (g^\circ,b^\circ) &&\text{in }C_{\mathrm{loc}}((0,\infty))^2,\\
            (g_j',b_j') &\rightharpoonup ((g^\circ)',(b^\circ)') &&\text{in }L^2(\delta,R)^2
      \end{align}
      for every \(0<\delta<R<\infty\).
      In particular, \(u^\circ \in H^1_{\mathrm{loc}}((0,\infty))^2 \subset AC_{\mathrm{loc}}(0,\infty)^2\).

      The local uniform convergence implies
      \begin{align}
            v \leq u^\circ \leq \mathbf1 \qquad\text{in }(0,\infty).
      \end{align}
      Passing to the limit in \(0\leq b_j(r)\leq C_0r\), we also obtain
      \begin{align}
            0 \leq b^\circ(r) \leq C_0r \qquad\text{for all }r>0,
      \end{align}
      and hence \(b^\circ(0^+)=0\).
      By weak lower semicontinuity of the derivative terms and local uniform convergence of the remaining terms, the energy of \(u^\circ\) on any compact subinterval of \((0,\infty)\) is bounded above by \(\liminf_{j\to\infty}E_\lambda(g_j,b_j)\).
      Since the energy density is nonnegative, letting the intervals increase to \((0,\infty)\) and applying the monotone convergence theorem, we obtain
      \begin{align}
            E_\lambda(g^\circ,b^\circ)
            &\leq \liminf_{j\to\infty}E_\lambda(g_j,b_j)\\
            &= \inf_{\mathcal K}E_\lambda < \infty.
      \end{align}
      Thus, \(u^\circ\in\mathcal K\) and \(u^\circ\) is a minimizer of \(E_\lambda\) over \(\mathcal K\).

      \textbf{Step 2.}
      In this step, we show that the minimizer \(u^\circ\) obtained in Step 1 satisfies the Euler--Lagrange equations \eqref{eq:EL}.
      Fix \(\phi:=(\phi_1,\phi_2)\in\Cc^2\) and \(t\in(0,1]\), and set
      \begin{align}
            z &:= u^\circ + t\phi=(z_1,z_2),\\
            y_0 &:= (|z_1|, \min\{z_2,1\}),\\
            y_i &:= y_{i-1}\vee v_i \qquad(1\leq i\leq m),\\
            \widehat z &:= (\min\{(y_m)_1,1\}, (y_m)_2).
      \end{align}
      All maxima and minima below are taken componentwise.
      Since \(z\) differs from \(u^\circ\) only on \(\supp\phi\), it has finite energy.      
      By Lemma~\ref{lem:lattice}(i), we have \(E_\lambda(y_0)\leq E_\lambda(z)<\infty\).
      Since each \(v_i\) has finite energy, repeated application of Lemma~\ref{lem:lattice}(ii) gives
      \begin{align}
            E_\lambda(y_i) \leq E_\lambda(y_{i-1}) + E_\lambda(v_i) < \infty
            \qquad(1\leq i\leq m).
      \end{align}
      Moreover, Lemma~\ref{lem:lattice}(i) implies \(E_\lambda(\widehat z) \leq E_\lambda(y_m) < \infty\).
      Put \(\Omega:=(0,\infty)\setminus\supp\phi\).
      Since \(z=u^\circ\) and \(0\leq u^\circ\leq\mathbf1\) in \(\Omega\), we have \(y_0=u^\circ\) in \(\Omega\).
      Moreover, since \(u^\circ\geq v_i\) for every \(i\), the definition of \(y_i\) implies \(y_i=u^\circ\) in \(\Omega\) for every \(i\).
      Consequently, \(\widehat z=u^\circ\) in \(\Omega\).
      By the definition of \(y_i\), we have
      \begin{align}
            y_m = y_0 \vee v_1 \vee \cdots \vee v_m
            = y_0 \vee v \geq v.
      \end{align}
      Since \(v\leq\mathbf1\), the first component of \(\widehat z\) satisfies
      \begin{align}
            (v)_1 &\leq \min\{(y_m)_1,1\} =(\widehat z)_1 \leq 1.
      \end{align}
      Moreover, since \((y_0)_2\leq1\) and \((v_i)_2\leq1\) for every \(i\), we have
      \begin{align}
            (v)_2 &\leq (y_m)_2 =(\widehat z)_2 \leq 1.
      \end{align}
      Consequently,
      \begin{align}
            v \leq \widehat z \leq \mathbf1 \qquad\text{in }(0,\infty).
      \end{align}
      Since \(\supp\phi\) is a compact subset of \((0,\infty)\) and \(\widehat z=u^\circ\) outside \(\supp\phi\), we have
      \begin{align}
            \lim_{r\downarrow0}(\widehat z)_2(r) = \lim_{r\downarrow0}b^\circ(r) = 0.
      \end{align}
      Thus, \(\widehat z\in\mathcal K\).

      For each \(i\in\{1,\ldots,m\}\), put
      \begin{align}\label{eq:theta}
            \vartheta^{(i)} := (v_i-y_{i-1})^+,
      \end{align}
      where the positive part is taken componentwise, that is,
      \begin{align}
            w^+ := (\max\{w_1,0\},\max\{w_2,0\})
      \end{align}
      for \(w=(w_1,w_2)\).
      Since \(y_{i-1} = u^\circ \geq v_i\) outside \(\supp\phi\), we have \(\supp\vartheta^{(i)}\subset\supp\phi\).
      Moreover, since \(v_i,y_{i-1} \in H^1_{\mathrm{loc}}((0,\infty))^2\), we obtain
      \begin{align}
            \vartheta^{(i)} \in (H_c^1)^2, \qquad
            \vartheta_1^{(i)} \geq 0
            \qquad\text{in }(0,\infty).
      \end{align}
      By the definition of the componentwise minimum,
      \begin{align}
            y_{i-1} \wedge v_i = v_i - \vartheta^{(i)}.
      \end{align}
      Since \(g^\circ\geq0\) and \(b^\circ\leq1\), we have
      \begin{align}
            |g^\circ + t\phi_1| &\geq g^\circ - t|\phi_1|,\\
            \min\{b^\circ + t\phi_2,1\} &\geq b^\circ - t|\phi_2|.
      \end{align}
      Thus, for \(k=1,2\), we have
      \begin{align}
            (y_0)_k \geq (u^\circ)_k - t|\phi_k|.
      \end{align}
      Since \(y_{i-1} \geq y_0\) and \(u^\circ \geq v_i\), it follows that
      \begin{align}\label{eq:thetabound}
            0 \leq \vartheta_k^{(i)} \leq t|\phi_k| \qquad(k=1,2)
      \end{align}
      in \((0,\infty)\).
      Since
      \begin{align}
            (y_{i-1})_1 \geq 0, \quad (v_i)_1 \geq 0, \qquad
            (y_{i-1})_2 \leq 1, \quad (v_i)_2 \leq 1
            \qquad\text{in }(0,\infty),
      \end{align}
      Lemma~\ref{lem:lattice}(ii) gives
      \begin{align}
            E_\lambda(y_i)
            &\leq E_\lambda(y_{i-1}) + E_\lambda(v_i) - E_\lambda(y_{i-1} \wedge v_i)\\
            &= E_\lambda(y_{i-1}) + E_\lambda(v_i) - E_\lambda(v_i - \vartheta^{(i)}).
            \label{eq:ys1}
      \end{align}
      We claim that there exists \(C_\phi>0\), independent of \(t\in(0,1]\) and \(i\in\{1,\ldots,m\}\), such that
      \begin{align}\label{eq:remainder}
            E_\lambda(v_i) - E_\lambda(v_i - \vartheta^{(i)})
            \leq \dd E_\lambda(v_i)[\vartheta^{(i)}] + C_\phi t^2.
      \end{align}
      To prove this inequality, write \(E_\lambda = T+S\), where
      \begin{align}
            T(g,b) &:= \frac{1}{2}\int_0^\infty \left\{|g'|^2 + \frac{n^2}{r^2}|b'|^2\right\}r\,\dd r,\\
            S(g,b) &:= \int_0^\infty \left\{\frac{n^2}{2r^2}g^2(1-b)^2 + \frac{\lambda}{8}(1-g^2)^2\right\}r\,\dd r.
      \end{align}
      Expanding the squares in \(T\), we obtain
      \begin{align}
            T(v_i) - T(v_i-\vartheta^{(i)})
            &=\dd T(v_i)[\vartheta^{(i)}]
            - \frac{1}{2}\int_0^\infty \left\{|(\vartheta_1^{(i)})'|^2 + \frac{n^2}{r^2}|(\vartheta_2^{(i)})'|^2\right\}r\,\dd r\\
            &\leq \dd T(v_i)[\vartheta^{(i)}].
      \end{align}
      
      Next, put \(K_\phi := \supp\phi \Subset (0,\infty)\).
      Taylor's formula gives
      \begin{align}
            S(v_i) - S(v_i-\vartheta^{(i)})
            &= \dd S(v_i)[\vartheta^{(i)}] - \int_0^1(1-s) D^2 S(v_i - s\vartheta^{(i)}) [\vartheta^{(i)}, \vartheta^{(i)}]\,\dd s.
      \end{align}
      By direct calculation,
      \begin{align}
            D^2S(g,b) [\vartheta^{(i)},\vartheta^{(i)}]
            &= \int_0^\infty \bigg\{
                  \left[
                        \frac{n^2}{r^2}(1-b)^2 + \frac{\lambda}{2}(3g^2-1)
                  \right](\vartheta_1^{(i)})^2\\
                  &\qquad
                  - \frac{4n^2}{r^2}g(1-b)\vartheta_1^{(i)}\vartheta_2^{(i)}
                  + \frac{n^2}{r^2}g^2 (\vartheta_2^{(i)})^2
            \bigg\}r\,\dd r.
      \end{align}
      Since \(0\leq v_i\leq\mathbf1\), \eqref{eq:thetabound} gives
      \begin{align}
            |(v_i)_k(r) - s\vartheta_k^{(i)}(r)|
            &\leq 1 + t|\phi_k(r)|\\
            &\leq 1 + \|\phi\|_{L^\infty}
            \qquad(k=1,2)
      \end{align}
      for \(r\in K_\phi\), \(s\in[0,1]\), and \(t\in(0,1]\).
      Since \(\supp\vartheta^{(i)} \subset K_\phi \Subset (0,\infty)\), the coefficients in the integrand above, evaluated at \((g,b) = v_i - s\vartheta^{(i)}\), are uniformly bounded for \(1\leq i\leq m\), \(s\in[0,1]\), and \(t\in(0,1]\).
      Hence, by \eqref{eq:thetabound},
      \begin{align}
            \left|D^2S(v_i-s\vartheta^{(i)})[\vartheta^{(i)},\vartheta^{(i)}]\right|
            &\leq C\int_{K_\phi}|\vartheta^{(i)}|^2r\,\dd r\\
            &\leq Ct^2\int_{K_\phi}|\phi|^2r\,\dd r
            \leq C_\phi t^2,
      \end{align}      
      where \(C,C_\phi>0\) are constants independent of \(i\), \(s\), and \(t\).
      Consequently,
      \begin{align}
            S(v_i) - S(v_i-\vartheta^{(i)})
            \leq \dd S(v_i)[\vartheta^{(i)}] + C_\phi t^2.
      \end{align}
      Combining this with the estimate for \(T\), we obtain \eqref{eq:remainder}.

      By Lemma~\ref{lem:sub}, the first variation satisfies
      \begin{align}
            \dd E_\lambda(v_i)[\vartheta^{(i)}] \leq 0,
      \end{align}
      because \(\vartheta_1^{(i)} \geq 0\).
      Applying \eqref{eq:remainder} to \eqref{eq:ys1}, we obtain
      \begin{align}
            E_\lambda(y_i)\leq E_\lambda(y_{i-1}) + C_\phi t^2 \qquad(1\leq i\leq m).
      \end{align}
      Iterating this inequality gives
      \begin{align}
            E_\lambda(y_m) \leq E_\lambda(y_0) + mC_\phi t^2.
      \end{align}
      Since \(\widehat z \in\mathcal K\) and \(u^\circ\) minimizes the energy over \(\mathcal K\), Lemma~\ref{lem:lattice}(i) yields
      \begin{align}
            E_\lambda(u^\circ)
            &\leq E_\lambda(\widehat z)\\
            &\leq E_\lambda(y_m)\\
            &\leq E_\lambda(y_0) + mC_\phi t^2\\
            &\leq E_\lambda(u^\circ + t\phi) + mC_\phi t^2.
      \end{align}
      Dividing by \(t>0\) and letting \(t\downarrow0\), we obtain
      \begin{align}
            \dd E_\lambda(u^\circ)[\phi] \geq 0.
      \end{align}
      Replacing \(\phi\) by \(-\phi\), we obtain the reverse inequality.
      Therefore,
      \begin{align}
            \dd E_\lambda(u^\circ)[\phi] = 0 \qquad\text{for all }\phi\in\Cc^2.
      \end{align}      
      This implies that \(u^\circ\) satisfies \eqref{eq:EL} in the weak sense on \((0,\infty)\).
      Since \(u^\circ \in H^1_{\mathrm{loc}}((0,\infty))^2\), standard regularity theory for ODEs yields \(u^\circ \in C^\infty((0,\infty))^2\).
     
      \textbf{Step 3.}
      In this step, we show that \(u^\circ\) satisfies the boundary conditions \eqref{eq:bc} and that its first component is strictly positive on \((0, \infty)\).
      We first verify that \(u^\circ\) satisfies \eqref{eq:bc}.
      Since \(v\leq u^\circ \leq \mathbf1\) and \(v(r) \to \mathbf1\) as \(r\to\infty\), we obtain
      \begin{align}
            \lim_{r\to\infty}(g^\circ(r),b^\circ(r)) = (1,1).
      \end{align}
      We have already shown that \(b^\circ(0^+) = 0\) and \(0 \leq b^\circ(r) \leq C_0r\).
      Hence we may choose \(\rho>0\) such that
      \begin{align}
            0 \leq b^\circ(r) \leq \frac{1}{2} \qquad\text{for all }r\in(0,\rho).
      \end{align}
      Using the terms involving \((g^\circ)'\) and \((g^\circ)^2(1-b^\circ)^2/r^2\) in \eqref{eq:radialenergy}, together with \(1-b^\circ\geq1/2\), the finite-energy bound gives
      \begin{align}\label{eq:ggprime}
            \int_0^\rho \left\{\frac{(g^\circ(r))^2}{r} + |(g^\circ)'(r)|^2r\right\}\dd r < \infty.
      \end{align}
      For \(0<r<s<\rho\), the Cauchy--Schwarz inequality implies
      \begin{align}
            |(g^\circ(s))^2 - (g^\circ(r))^2|
            &\leq 2\int_r^s |g^\circ(\tau)(g^\circ)'(\tau)|\,\dd\tau\\
            &\leq 2\left(\int_r^s\frac{(g^\circ(\tau))^2}{\tau}\,\dd\tau\right)^{1/2}
                   \left(\int_r^s|(g^\circ)'(\tau)|^2\tau\,\dd\tau\right)^{1/2}.
      \end{align}
      By \eqref{eq:ggprime}, the right-hand side tends to zero as \(r,s\downarrow0\).
      Therefore, \((g^\circ(r))^2\) has a finite limit as \(r\downarrow0\).
      If this limit were positive, then \(\int_0^\rho (g^\circ(r))^2/r\,\dd r\) would diverge.
      This contradicts the finiteness of the integral established in \eqref{eq:ggprime}.
      Thus,
      \begin{align}
            \lim_{r\downarrow0}(g^\circ(r))^2 = 0.
      \end{align}
      Since \(g^\circ \geq 0\), it follows that \(g^\circ(0^+) = 0\).
      We have shown that \(u^\circ \in C^2((0,\infty))^2\) satisfies \eqref{eq:EL} and \eqref{eq:bc}, so it is a radial solution.
      Finally, \(g^\circ\geq0\) and Lemma~\ref{lem:shape}(iv) imply
      \begin{align}
            0 < g^\circ < 1 \qquad\text{in }(0,\infty).
      \end{align}
      Hence, \(u^\circ\) is a positive radial solution satisfying \(u^\circ \geq v\), as asserted.
\end{proof}

\begin{lemma}\label{lem:compare}
      Let \(u^\circ = (g^\circ, b^\circ)\) be a positive radial solution, and let \((f,a)\) be any radial solution such that
      \begin{align}
            |f| \leq g^\circ, \qquad a \leq b^\circ \quad \text{in }(0,\infty).
      \end{align}
      Then \((f,a) = u^\circ\).
\end{lemma}

\begin{proof}
      Let
      \begin{align}
            B := \frac{na'}{r}, \qquad
            B^\circ := \frac{n(b^\circ)'}{r},
      \end{align}
      and define
      \begin{align}
            \Theta := n(b^\circ - a).
      \end{align}
      By the assumption \(a\leq b^\circ\), we have \(\Theta\geq0\) in \((0,\infty)\).
      Moreover, the boundary conditions \eqref{eq:bc} give
      \begin{align}
            \Theta(0^+) = 0, \qquad \lim_{r\to\infty}\Theta(r) = 0.
      \end{align}
      Since \(0<a,b^\circ<1\) by Lemma~\ref{lem:shape}(i), we also have
      \begin{align}
            0 \leq \Theta(r) \leq n \qquad\text{for all }r \in (0,\infty).
      \end{align}
      Furthermore,
      \begin{align}
            \Theta'(r) = r(B^\circ(r) - B(r)) \qquad\text{for all }r \in (0,\infty).
      \end{align}

      We first compare \(B\) and \(B^\circ\).
      By Lemma~\ref{lem:ends}, we have \(B, B^\circ \in L^2((0,\infty),r\dd r)\).
      Moreover, Lemmas~\ref{lem:shape}(i) and \ref{lem:ends}(i) imply that \(B\) is decreasing and
      \begin{align}
            B(0^+) \in (0,\infty), \qquad
            \lim_{r\to\infty}B(r) = 0.
      \end{align}
      Hence, \(B' \in L^1(0,\infty)\).
      Since \(\Theta(0^+) = \Theta(\infty) = 0\), integration by parts over \((\varepsilon,R)\), followed by \(\varepsilon\downarrow0\) and \(R\to\infty\), gives
      \begin{align}\label{eq:compareint}
            2\int_0^\infty B(B^\circ - B)r\dd r
            &= 2\int_0^\infty B\Theta'\dd r\\
            &= -2\int_0^\infty B'\Theta\dd r
            \geq 0,
      \end{align}
      where the last inequality follows from \(\Theta\geq0\) and the fact that \(B'\leq0\) by Lemma~\ref{lem:shape}(i).
      Using the identity
      \begin{align}
            (B^\circ)^2 - B^2 = (B^\circ - B)^2 + 2B(B^\circ - B),
      \end{align}
      and integrating it with respect to \(r\dd r\), we obtain from \eqref{eq:compareint}
      \begin{align}\label{eq:compareB}
            \int_0^\infty (B^\circ)^2r\,\dd r
            \geq \int_0^\infty B^2r\,\dd r + \int_0^\infty (B^\circ-B)^2r\,\dd r.
      \end{align}
      Next, we compare the potential terms.
      Since \(u^\circ\) is a positive radial solution, Lemma~\ref{lem:shape}(iv) gives
      \begin{align}
            0 < g^\circ < 1 \qquad\text{in }(0,\infty).
      \end{align}
      Together with the assumption \(|f| \leq g^\circ\), this implies
      \begin{align}
            (1-f(r)^2)^2 \geq (1-(g^\circ(r))^2)^2 \qquad\text{for all }r \in (0,\infty).
      \end{align}
      Applying the identity \eqref{eq:virial} to \((f,a)\) and to \(u^\circ\), we obtain
      \begin{align}
            \int_0^\infty B^2r\,\dd r
            &= \frac{\lambda}{4}\int_0^\infty (1-f^2)^2r\,\dd r\\
            &\geq \frac{\lambda}{4}\int_0^\infty (1-(g^\circ)^2)^2r\,\dd r\\
            &= \int_0^\infty (B^\circ)^2r\,\dd r.
      \end{align}
      Combining this inequality with \eqref{eq:compareB}, we find
      \begin{align}
            \int_0^\infty B^2r\,\dd r
            &\geq \int_0^\infty (B^\circ)^2r\,\dd r\\
            &\geq \int_0^\infty B^2r\,\dd r + \int_0^\infty (B^\circ-B)^2r\,\dd r.
      \end{align}
      Therefore,
      \begin{align}
            \int_0^\infty (B^\circ-B)^2r\,\dd r = 0.
      \end{align}
      Since \(B\) and \(B^\circ\) are continuous, it follows that \(B=B^\circ\) in \((0,\infty)\).
      Applying \eqref{eq:virial} to both solutions, we obtain
      \begin{align}
            &\int_0^\infty \left\{(1-f^2)^2 - (1-(g^\circ)^2)^2\right\}r\,\dd r\\
            &= \frac{4}{\lambda} \int_0^\infty \left\{B^2 - (B^\circ)^2\right\}r\,\dd r
            = 0.
      \end{align}
      Since \(|f| \leq g^\circ < 1\), the integrand on the left is nonnegative and vanishes only when \(|f| = g^\circ\).
      Hence, \(|f| = g^\circ > 0\) in \((0,\infty)\).
      Since \(f\) is continuous and \(f(r)\to1\) as \(r\to\infty\), it follows that
      \begin{align}
            f = g^\circ \qquad\text{in }(0,\infty).
      \end{align}

      Finally, since \(B=B^\circ\), we have \(a' = (b^\circ)'\) in \((0,\infty)\).
      Together with the boundary conditions \(a(0^+) = b^\circ(0^+)=0\), this yields
      \begin{align}
            a = b^\circ \qquad\text{in }(0,\infty).
      \end{align}
      Thus, \((f,a) = (g^\circ, b^\circ) = u^\circ\), as asserted.
\end{proof}

\begin{proof}[Proof of Theorem~\ref{thm:A}]
      Existence follows from \cite{MR987294}.
      We prove uniqueness.
      Let \(u_1=(f_1,a_1)\) and \(u_2=(f_2,a_2)\) be any two radial solutions.
      By Lemma~\ref{lem:upper}, there exists a positive radial solution \(u^\circ = (g^\circ,b^\circ)\) such that
      \begin{align}
            g^\circ \geq \max\{|f_1|,|f_2|\},\qquad
            b^\circ \geq \max\{a_1,a_2\}
            \quad\text{in }(0,\infty).
      \end{align}
      Thus, the assumptions of Lemma~\ref{lem:compare} are satisfied by both \(u_1\) and \(u_2\).
      Applying that lemma to each solution, we obtain
      \begin{align}
            u_1 = u^\circ = u_2,
      \end{align}
      which proves uniqueness.
\end{proof}

\section{Analyticity and monotonicity with respect to the coupling constant}
As a consequence of Theorem~\ref{thm:A}, we show that the unique radial solution depends real analytically on the coupling constant \(\lambda\) and that both components are strictly increasing in \(\lambda\).
The key ingredient is the positivity of the radial Hessian.
Fix \(n\geq1\) and \(\lambda>0\), and let \((f,a)\) denote the unique radial solution given by Theorem~\ref{thm:A}.
We first introduce the radial Hessian of \(E_\lambda\) at \((f,a)\).
For \(\xi,\eta\in\Cc\), consider variations of the form
\begin{align}
      \delta f := \xi, \qquad \delta a:= \frac{r}{n}\eta.
\end{align}
A direct calculation of the second variation of \eqref{eq:radialenergy}, using integration by parts, gives
\begin{align}\label{eq:Q}
      Q[\xi,\eta]
      &:= \left.\frac{\dd^2}{\dd t^2} E_\lambda\left(f + t\xi, a + t\frac{r}{n}\eta\right)\right|_{t=0}\\
      &= \int_0^\infty \left\{
            |\xi'|^2 + V_{11}\xi^2 + |\eta'|^2 + V_{22}\eta^2 + 2V_{12}\xi\eta
      \right\}r\dd r,
\end{align}
where the entries of the symmetric matrix potential \(V\) are given by
\begin{align}\label{eq:V}
      V_{11} &:= n^2P^2 + \frac{\lambda}{2}(3f^2-1),\\
      V_{22} &:= \frac{1}{r^2} + f^2,\\
      V_{12} &= V_{21} := -2nPf.
\end{align}
We denote by \(M_0\) the corresponding matrix differential expression
\begin{align}\label{eq:M0}
      M_0 &:= -\Del I_2 + V
      = \begin{pmatrix}
            -\Del + V_{11} & V_{12}\\
            V_{21} & -\Del + V_{22}
      \end{pmatrix}.
\end{align}
For compactly supported smooth variations, integration by parts yields
\begin{align}
      Q[\xi,\eta] = \int_0^\infty
      \begin{pmatrix}\xi & \eta\end{pmatrix}
      M_0 \begin{pmatrix}\xi \\ \eta\end{pmatrix}
      r\dd r.
\end{align}
We next establish an identity for \(M_0\) that will be used to prove positivity of the radial Hessian.

\begin{lemma}\label{lem:D}
      Define
      \begin{align}
            D := (D_1,D_2) := (rf',na').
      \end{align}
      Then the following properties hold:
      \begin{align}
            D_1(r) > 0, \qquad D_2(r) > 0
            \qquad\text{for all }r\in(0,\infty).
      \end{align}
      Moreover, the following identity holds pointwise in \((0,\infty)\):
      \begin{align}\label{eq:MD}
            M_0 D
            = \begin{pmatrix}
                  \lambda f(1-f^2)\\
                  2nPf^2
            \end{pmatrix}.
      \end{align}
      Both components on the right-hand side are strictly positive in \((0,\infty)\).
\end{lemma}
\begin{proof}
      The positivity of \(D\) follows from Lemma~\ref{lem:shape}.
      A direct calculation using \eqref{eq:EL} gives the stated identity for \(M_0D\).
      Its right-hand side is positive since \(0<f<1\) and \(P>0\) by Lemma~\ref{lem:shape}.
\end{proof}

We next derive a representation of the quadratic form \(Q\) using the positive vector \(D\) from Lemma~\ref{lem:D}.

\begin{lemma}\label{lem:gst}
      For any \((\xi,\eta)\in\Cc^2\), put
      \begin{align}
            p := \frac{\xi}{D_1}, \qquad q := \frac{\eta}{D_2}.
      \end{align}
      Then \(p,q\in\Cc\), and the following identity holds:
      \begin{align}\label{eq:gst}
            Q[\xi,\eta] = \int_0^\infty \bigl\{
                  &D_1^2(p')^2 + D_2^2(q')^2 + 2nPfD_1D_2(p-q)^2\\
                  &+ \lambda f(1-f^2)D_1p^2 + 2nPf^2D_2q^2\bigr\}r\dd r.
      \end{align}
\end{lemma}
\begin{proof}
      Since \(D_1,D_2>0\) by Lemma~\ref{lem:D}, we have \(p,q\in\Cc\).
      For \(j=1,2\) and \(h\in\Cc\), integration by parts gives
      \begin{align}
            \int_0^\infty ((D_jh)')^2r\dd r 
            = \int_0^\infty \left\{D_j^2(h')^2 - D_j(\Del D_j)h^2\right\}r\dd r.
      \end{align}
      Substituting \(\xi=D_1p\) and \(\eta=D_2q\) into \eqref{eq:Q},
      we therefore obtain
      \begin{align}
            Q[\xi,\eta]
            =\int_0^\infty\bigl\{
                  &D_1^2(p')^2 + D_2^2(q')^2 + (M_0 D)_1 D_1 p^2 + (M_0 D)_2 D_2 q^2\\
                  &- V_{12}D_1 D_2 (p-q)^2
            \bigr\}r\dd r.
      \end{align}
      Using \eqref{eq:MD} and \(V_{12}=-2nPf\), we obtain \eqref{eq:gst}.
\end{proof}

We now establish strict positivity of the radial Hessian.
Let
\begin{align}
      X := L^2((0,\infty),r\dd r),
\end{align}
and define
\begin{align}
      \mathcal H&:=\overline{\Cc}^{\|\cdot\|_{H_1}},\\
      \|h\|_{H_m}^2
      &:=\int_0^\infty\left\{
            |h'|^2+\frac{m^2}{r^2}|h|^2+|h|^2
      \right\}r\dd r
      \qquad(m\geq1).
\end{align}
For each integer \(m\geq1\), let \(H_m\) denote the space \(\mathcal H\) equipped with the norm \(\|\cdot\|_{H_m}\).
These norms are equivalent, since
\begin{align}
      \|h\|_{H_1}\leq\|h\|_{H_m}\leq m\|h\|_{H_1}.
\end{align}
We also set \(H:=H_n\times H_1\), equipped with the product norm.

\begin{proposition}\label{prop:hessian}
      The quadratic form \(Q\) defined in \eqref{eq:Q} on \(\Cc^2\) is lower semibounded and closable in \(X\oplus X\), and its closure has form domain \(H\).
      Let \(M_0\) denote the self-adjoint operator associated with this closed form.
      Then there exists \(c>0\) such that
      \begin{align}\label{eq:hessianpositive}
            M_0 \geq cI \qquad\text{on }X\oplus X.
      \end{align}
      In particular, \(\ker M_0=\{0\}\).
\end{proposition}
\begin{proof}
      We divide the proof into three steps.
      We first identify the domain of the closed quadratic form.
      We then prove nonnegativity and the triviality of the kernel.
      Finally, we establish the positive spectral gap.

      \textbf{Step 1.}
      We first show that \(Q\) is lower semibounded and closable in \(X\oplus X\), and that its closure has form domain \(H\).
      By \eqref{eq:V}, the matrix potential can be decomposed as
      \begin{align}
            V(r) = \begin{pmatrix}
                  \frac{n^2}{r^2} & 0\\
                  0 & \frac{1}{r^2}
            \end{pmatrix} + \widetilde V(r),
      \end{align}
      where
      \begin{align}
            \widetilde V(r) := \begin{pmatrix}
                  \frac{n^2(a^2-2a)}{r^2} + \frac{\lambda}{2}(3f^2-1) & -2nPf\\
                  -2nPf & f^2
            \end{pmatrix}.
      \end{align}
      We claim that \(\widetilde V\) is bounded in \((0,\infty)\).
      Indeed, Lemma~\ref{lem:ends}(i) gives \(a(r)=O(r^2)\) and \(f(r)=O(r^n)\) as \(r\downarrow0\).
      Since \(P=(1-a)/r\) and \(n\geq1\), these estimates imply that both \((a^2-2a)/r^2\) and \(Pf\) are bounded near the origin.
      On the other hand, Lemma~\ref{lem:shape} gives \(0<f,a<1\) in \((0,\infty)\), so the entries of \(\widetilde V\) are also bounded away from the origin.
      Thus, we may set
      \begin{align}
            M := \sup_{r>0}\|\widetilde V(r)\|_{\mathrm{op}} < \infty.
      \end{align}
      For \(k := (\xi,\eta)\in\Cc^2\), the decomposition above gives
      \begin{align}
            Q[k] = \|k\|_H^2-\|k\|_{X\oplus X}^2 + \int_0^\infty k^{\mathsf T}\widetilde V k\,r\dd r.
      \end{align}
      Since \(\widetilde V\) is symmetric and bounded, we have
      \begin{align}
            \left|\int_0^\infty k^{\mathsf T}\widetilde V k\,r\dd r\right|
            \leq M\|k\|_{X\oplus X}^2.
      \end{align}
      Consequently,
      \begin{align}\label{eq:garding}
            \|k\|_H^2
            &\leq Q[k] + (M+1)\|k\|_{X\oplus X}^2\\
            &\leq (1 + 2M)\|k\|_H^2
            \qquad\text{for all }k\in\Cc^2.
      \end{align}
      In particular, \(Q\) is lower semibounded.
      Moreover, the norm
      \begin{align}
            \|k\|_Q := \bigl(Q[k] + (M+1)\|k\|_{X\oplus X}^2\bigr)^{1/2}
      \end{align}
      is equivalent to \(\|k\|_H\) on \(\Cc^2\).
      Since \(\Cc^2\) is dense in \(H\), the two norms are equivalent, and \(H\) is continuously embedded in \(X \oplus X\), the form \(Q\) is closable in \(X \oplus X\), and its closure has form domain \(H\).
      We also use \(Q\) to denote the closed quadratic form.

      \textbf{Step 2.}
      We next show that \(Q\) is nonnegative on \(H\) and that the kernel of \(M_0\) is trivial.
      By Lemma~\ref{lem:gst}, every term on the right-hand side of \eqref{eq:gst} is nonnegative.
      Hence \(Q[k]\geq0\) for all \(k\in\Cc^2\), and this inequality extends to all \(k\in H\) by density.
      It follows that \(M_0\geq0\).

      We next show that \(\ker M_0=\{0\}\).
      Suppose that \(k := (\xi,\eta)\in\ker M_0\).
      Since \(k\in H\), there exists a sequence \(k_j:=(\xi_j,\eta_j)\in\Cc^2\) converging to \(k\) in \(H\).
      By the continuity of the closed form on \(H\), we have
      \begin{align}
            Q[k_j] \longrightarrow Q[k]
            = \langle M_0k, k\rangle_{X \oplus X}=0.
      \end{align}
      Fix \(0<\delta<R<\infty\).
      By \eqref{eq:gst}, we have
      \begin{align}
            Q[k_j] &\geq \int_\delta^R \left\{
                  \lambda f(1-f^2)D_1p_j^2 + 2nPf^2D_2q_j^2
            \right\}r\dd r\\
            &= \int_\delta^R \left\{
                  \frac{\lambda f(1-f^2)}{D_1}|\xi_j|^2 + \frac{2nPf^2}{D_2}|\eta_j|^2
            \right\}r\dd r.
      \end{align}
      By Lemma~\ref{lem:D}, the coefficients in the last integral are continuous and strictly positive on \([\delta,R]\).
      Hence, they are bounded below by a positive constant \(c_{\delta,R}>0\), and we obtain
      \begin{align}
            Q[k_j] \geq c_{\delta,R}\int_\delta^R \bigl(|\xi_j(r)|^2 + |\eta_j(r)|^2\bigr)r\dd r.
      \end{align}
      Letting \(j\to\infty\), we obtain
      \begin{align}
            \int_\delta^R \bigl(|\xi(r)|^2 + |\eta(r)|^2\bigr)r\dd r=0.
      \end{align}      
      Since \(\delta\) and \(R\) are arbitrary, we conclude that \(k=0\).
      Thus, \(\ker M_0=\{0\}\).
      
      \textbf{Step 3.}
      Finally, we show that the spectrum of \(M_0\) is bounded away from zero.
      The boundary conditions \eqref{eq:bc} imply
      \begin{align}
            \widetilde V(r) \longrightarrow
            \begin{pmatrix}
                  \lambda & 0\\
                  0 & 1
            \end{pmatrix}
            \qquad\text{as }r\to\infty.
      \end{align}
      Let \(A_0\) be the self-adjoint operator in \(X\oplus X\) associated with the quadratic form
      \begin{align}
            Q_0[\xi,\eta] := \int_0^\infty \left\{
                  |\xi'|^2 + \left(\frac{n^2}{r^2} + \lambda\right)|\xi|^2 + |\eta'|^2 + \left(\frac1{r^2} + 1\right)|\eta|^2
            \right\}r\dd r
      \end{align}
      with form domain \(H\).
      The standard spectral description of the radial Laplacian, for example by the Hankel transform, gives
      \begin{align}
            \sigma(A_0) = \sigma_{\mathrm{ess}}(A_0) = [\min\{\lambda,1\},\infty).
      \end{align}

      Put
      \begin{align}
            W(r) := \widetilde V(r) - \begin{pmatrix} \lambda&0 \\ 0&1\end{pmatrix}.
      \end{align}
      Then \(W\) is bounded and \(W(r)\to0\) as \(r\to\infty\), and the associated operators satisfy \(M_0 = A_0 + W\).
      We claim that multiplication by \(W\) is compact from \(H\) into \(X\oplus X\).
      To see this, first note that for any \(h\in H_m\), the definition of the norm gives
      \begin{align}
            \int_0^\delta |h(r)|^2r\dd r
            \leq \frac{\delta^2}{m^2}\|h\|_{H_m}^2
            \qquad\text{for all }\delta>0.
      \end{align}
      Since \(W\) is bounded, the preceding estimate shows that \(Wk\) is uniformly small near the origin in \(X\oplus X\) for \(k\) in bounded subsets of \(H\).
      The same holds near infinity since \(W(r)\to0\) as \(r\to\infty\).
      On each interval \([\delta,R]\), the embedding \(H^1(\delta,R)\hookrightarrow L^2(\delta,R)\) is compact.
      Consequently, multiplication by \(W\) is compact from \(H\) into \(X\oplus X\).
      Since \(\lambda>0\), the spectral identity above implies that \(-1\) belongs to the resolvent set of \(A_0\).
      Moreover, the form norm associated with \(Q_0\) is equivalent to the norm of \(H\), so \((A_0+I)^{-1}\) maps \(X\oplus X\) boundedly into \(H\).
      Consequently, \(W(A_0+I)^{-1}\) is compact.
      Thus, \(W\) is a relatively compact perturbation of \(A_0\), and Weyl's theorem yields
      \begin{align}
            \sigma_{\mathrm{ess}}(M_0)
            = \sigma_{\mathrm{ess}}(A_0)
            = [\min\{\lambda,1\},\infty).
      \end{align}

      Since \(\lambda>0\), the essential spectrum of \(M_0\) is bounded away from zero.
      Together with \(M_0\geq0\) and \(\ker M_0=\{0\}\) from Step 2,
      this implies that \(0\notin\sigma(M_0)\).
      Therefore, there exists \(c>0\) such that
      \begin{align}
            \inf\sigma(M_0)\geq c > 0.
      \end{align}
      This proves \eqref{eq:hessianpositive} and completes the proof.
\end{proof}

We now use Proposition~\ref{prop:hessian} to study the dependence of the radial solution on the coupling constant.
For each \(\lambda>0\), write \(u_\lambda := (f_\lambda,a_\lambda)\) for the unique radial solution given by Theorem~\ref{thm:A}.

\begin{proposition}\label{prop:analytic}
      Define
      \begin{align}
            f_*(r) := \left(\frac{r^2}{1+r^2}\right)^{n/2},\qquad
            a_*(r) := \frac{r^2}{1+r^2},
      \end{align}
      and, for \(k := (\xi,\eta)\in H\), put
      \begin{align}\label{eq:Phi}
            \Phi(k) := \left(f_* + \xi, a_* + \frac{r}{n}\eta\right).
      \end{align}
      Then the map
      \begin{align}\label{eq:analytic}
            (0,\infty)\ni\lambda\longmapsto
            \Phi^{-1}(u_\lambda)
            = \left(f_\lambda-f_*, \frac{n}{r}(a_\lambda-a_*)\right)\in H
      \end{align}
      is real analytic.
      Moreover, for every \(\lambda>0\) and \(r>0\), we have
      \begin{align}\label{eq:monotone}
            \partial_\lambda f_\lambda(r) > 0, \qquad
            \partial_\lambda a_\lambda(r) > 0.
      \end{align}
\end{proposition}
\begin{proof}
      We divide the proof into three steps.
      We first formulate the Euler--Lagrange equations as an analytic equation on \(H\).
      We then use Proposition~\ref{prop:hessian} and the analytic implicit-function theorem to establish the real analyticity of the radial solution with respect to \(\lambda\).
      Finally, we prove its strict monotonicity with respect to \(\lambda\).

      \textbf{Step 1.}
      We show that the first variation of the energy defines a real-analytic map from \((0,\infty)\times H\) to \(H^*\).
      For \(k\in H\), set
      \begin{align}
            \mathscr E(\lambda,k) &:= E_\lambda(\Phi(k)), &
            \mathscr F(\lambda,k) &:= \dd_k\mathscr E(\lambda,k)\in H^*.
      \end{align}
      Write \(k=(\xi,\eta)\) and put
      \begin{align}
            P_* := \frac{1-a_*}{r}, \qquad B_* := \frac{na_*'}{r}.
      \end{align}
      Substitution into \eqref{eq:radialenergy} gives
      \begin{align}
            \mathscr E(\lambda,k)
            &= \frac{1}{2}\|f_*' + \xi'\|_X^2 + \frac{1}{2}\|nP_*f_* + nP_*\xi - f_*\eta - \xi\eta\|_X^2\\
            &\quad +\frac{1}{2}\left\|B_* + \eta' + \frac\eta r\right\|_X^2 + \frac{\lambda}{8}\|1 - f_*^2 - 2f_*\xi - \xi^2\|_X^2.
            \label{eq:analyticenergy}
      \end{align}
      By the definitions of \(f_*\) and \(a_*\), we have
      \begin{align}
            f_*', \quad nP_*f_*, \quad B_*, \quad 1-f_*^2 \in X,
      \end{align}
      and \(0\leq f_*\leq1\), \(|nP_*|\leq n/r\).
      For any \(\xi,\widetilde{\xi}\in H_n\) and \(\eta\in H_1\), the Sobolev inequality and H\"older's inequality give
      \begin{align}
            \|\xi\eta\|_X
            &\leq \|\xi\|_{L^4}\|\eta\|_{L^4}
            \leq C\|\xi\|_{H_n}\|\eta\|_{H_1},\\
            \|\xi\widetilde{\xi}\|_X
            &\leq \|\xi\|_{L^4}\|\widetilde{\xi}\|_{L^4}
            \leq C\|\xi\|_{H_n}\|\widetilde{\xi}\|_{H_n},
      \end{align}
      where the \(L^4\)-norms are taken with respect to \(r\dd r\).
      Thus \(\mathscr E\) is a continuous polynomial in \((\lambda,k)\).
      In particular,
      \begin{align}
            \mathscr F : (0,\infty) \times H \longrightarrow H^*
      \end{align}
      is real analytic.

      By Lemma~\ref{lem:ends} and the explicit formulas for \(f_*\) and \(a_*\), we have
      \begin{align}
            k_\lambda := \Phi^{-1}(u_\lambda) \in H
            \qquad\text{for every }\lambda > 0.
      \end{align}
      Since \(u_\lambda\) satisfies \eqref{eq:EL}, its first variation vanishes on \(\Cc^2\).
      Since \(\Cc^2\) is dense in \(H\) and \(\mathscr F(\lambda,k_\lambda)\in H^*\), this identity extends to all of \(H\).
      Hence
      \begin{align}\label{eq:stationary}
            \mathscr F(\lambda, k_\lambda) = 0.
      \end{align}

      \textbf{Step 2.}
      We prove the real-analytic dependence of \(k_\lambda\) on \(\lambda\).
      Fix \(\lambda_0>0\).
      The derivative
      \begin{align}
            \dd_k\mathscr F(\lambda_0,k_{\lambda_0}) : H \longrightarrow H^*
      \end{align}
      is the bounded operator associated with the bilinear form of the radial Hessian \(Q\) at \(u_{\lambda_0}\).
      Proposition~\ref{prop:hessian} gives \(c>0\) such that
      \begin{align}
            Q[h]\geq c\|h\|_{X\oplus X}^2
            \qquad\text{for every }h\in H.
      \end{align}
      Combining this estimate with \eqref{eq:garding}, we obtain
      \begin{align}
            \|h\|_H^2 \leq Q[h] + (M+1)\|h\|_{X\oplus X}^2
            \leq \left(1 + \frac{M+1}{c}\right)Q[h].
      \end{align}
      Hence the bilinear form is coercive on \(H\), and the Lax--Milgram theorem implies that
      \(\dd_k\mathscr F(\lambda_0,k_{\lambda_0})\) is an isomorphism from \(H\) onto \(H^*\).
      The analytic implicit-function theorem therefore yields an interval \(I\) containing \(\lambda_0\) and a real-analytic map \(k:I\to H\) such that
      \begin{align}
            k(\lambda_0) = k_{\lambda_0},\qquad
            \mathscr F(\lambda, k(\lambda)) = 0
            \qquad\text{for all }\lambda \in I.
      \end{align}

      We show that this local branch coincides with the radial solutions.
      Fix \(\lambda\in I\), write \(k(\lambda) = (\xi,\eta) \in H\), and set \((f,a) := \Phi(k(\lambda))\).
      The equation \(\mathscr F(\lambda, k(\lambda)) = 0\) implies that \((f,a)\) satisfies \eqref{eq:EL} weakly in \((0,\infty)\).
      Standard regularity theory for ODEs gives \((f,a)\in C^\infty((0,\infty))^2\).
      Setting \(s:=\log r\), we have \(\xi(e^s),\eta(e^s)\in H^1(\R)\).
      Hence the Sobolev embedding \(H^1(\R)\hookrightarrow C_0(\R)\) implies
      \begin{align}
            \xi(r),\eta(r) \longrightarrow 0
            \qquad\text{as }r\downarrow0\text{ or }r\to\infty.
      \end{align}
      Therefore, \eqref{eq:Phi} gives
      \begin{align}
            f(0^+) = 0, \qquad f(\infty) = 1, \qquad a(0^+) = 0.
      \end{align}
      It remains to prove that \(a(r)\to1\) as \(r\to\infty\).

      Put \(w:=1-a\) and \(W:=ww'/r\).
      The second equation of \eqref{eq:EL} gives
      \begin{align}\label{eq:Wprime}
            W' = \frac{w'^2}{r} + \frac{f^2w^2}{r} \geq 0 \quad \text{in } (0,\infty).
      \end{align}
      Since \(k(\lambda)\in H\), \eqref{eq:analyticenergy} gives \(E_\lambda(f,a)<\infty\).
      Together with \(f(r)\to1\), this implies that there exists \(R>0\) such that
      \begin{align}\label{eq:wintegrable}
            \int_R^\infty \frac{w(r)^2}{r}\dd r < \infty.
      \end{align}
      If \(W(r_0)>0\) at some \(r_0>0\), then \eqref{eq:Wprime} implies
      \begin{align}
            (w^2)'(r) = 2rW(r) \geq 2rW(r_0)
            \qquad\text{for every }r \in [r_0,\infty).
      \end{align}
      Integrating from \(r_0\) to \(r\), we obtain
      \begin{align}
            w(r)^2 \geq W(r_0)(r^2-r_0^2)
            \qquad\text{for every }r \in [r_0,\infty).
      \end{align}
      Hence \(\int_R^\infty w(r)^2/r\,\dd r = \infty\), which contradicts \eqref{eq:wintegrable}.
      Therefore \(W\leq0\) in \((0,\infty)\), and \(w^2\) is nonincreasing.
      If \(w(r_1)=0\) for some \(r_1>0\), then the nonnegativity and monotonicity of \(w^2\) imply
      \begin{align}
            w(r) = 0 \qquad\text{for all }r \in [r_1,\infty).
      \end{align}
      In particular, \(w(r_1) = w'(r_1) = 0\).
      Since \(w\) satisfies
      \begin{align}
            w'' - \frac{w'}{r} = f^2w \qquad\text{in }(0,\infty),
      \end{align}
      the unique solvability of initial-value problems of ODE implies \(w\equiv0\) in \((0,\infty)\), which contradicts \(w(0^+)=1\).
      Thus \(w>0\) in \((0,\infty)\).
      Since \(w^2\) is nonnegative and nonincreasing, it has a limit as \(r\to\infty\).
      By \eqref{eq:wintegrable}, this limit must be zero.
      Since \(w=1-a\), we conclude that
      \begin{align}
            \lim_{r\to\infty}a(r) = 1.
      \end{align}
      Hence \((f,a)\) satisfies \eqref{eq:bc} and is a radial solution.
      By Theorem~\ref{thm:A}, we conclude that
      \begin{align}
            k(\lambda) = k_\lambda \qquad\text{for all }\lambda\in I.
      \end{align}

      Since \(\lambda_0>0\) was arbitrary, the map \(\lambda \mapsto k_\lambda = \Phi^{-1}(u_\lambda)\) is real analytic from \((0,\infty)\) into \(H\).

      \textbf{Step 3.}
      Finally, we establish the strict monotonicity of both profiles.
      Fix \(\lambda>0\) and write \(f:=f_\lambda\), \(a:=a_\lambda\).
      Define
      \begin{align}
            \zeta := \partial_\lambda k_\lambda
            = \left(\partial_\lambda f, \frac{n}{r}\partial_\lambda a\right)\in H.
      \end{align}
      Differentiating \eqref{eq:stationary} with respect to \(\lambda\), we obtain
      \(\dd_k\mathscr F(\lambda,k_\lambda)\zeta = -\partial_\lambda\mathscr F(\lambda,k_\lambda)\) in \(H^*\).
      A direct calculation of the right-hand side gives
      \begin{align}\label{eq:zeta}
            M_0\zeta = \begin{pmatrix}
                  \frac{1}{2}f(1-f^2)\\
                  0
            \end{pmatrix}
      \end{align}
      in the weak sense.
      Since \(f\) and \(a\) are smooth in \((0,\infty)\), the matrix potential \(V\) defined in \eqref{eq:M0} and the right-hand side of \eqref{eq:zeta} are smooth in \((0,\infty)\).
      Hence standard regularity theory for linear ODE systems yields \(\zeta\in C^\infty((0,\infty))^2\), and \eqref{eq:zeta} holds pointwise.

      Let \(\zeta^+\) denote the componentwise positive part of \(\zeta\), and set \(\zeta^- := (-\zeta)^+\).
      For \(k=(\xi,\eta)\in H\) and \(h=(\widetilde{\xi},\widetilde{\eta})\in H\), define
      \begin{align}
            Q(k,h) := \int_0^\infty \bigl\{
                  \xi'\widetilde{\xi}'
                  + \eta'\widetilde{\eta}'
                  + V_{11}\xi\widetilde{\xi}
                  + V_{22}\eta\widetilde{\eta}
                  + V_{12}(\xi\widetilde{\eta}
                  + \widetilde{\xi}\eta)
            \bigr\}r\dd r.
      \end{align}
      Note that \(Q(k,k)=Q[k]\) for every \(k\in H\).
      Testing \eqref{eq:zeta} against \(\zeta^-\), we obtain
      \begin{align}
            0 &\leq \int_0^\infty \frac{1}{2}f(1-f^2)\zeta_1^-r\dd r\\
            &= Q(\zeta,\zeta^-)
            = Q(\zeta^+,\zeta^-) - Q[\zeta^-]\\
            &\leq -Q[\zeta^-]
            \leq -c\|\zeta^-\|_{X\oplus X}^2.
      \end{align}
      Here the first inequality follows from \(0<f<1\), while \(Q(\zeta^+,\zeta^-)\leq 0\) follows from \(V_{12} = -2nPf<0\).
      The last inequality follows from Proposition~\ref{prop:hessian}.
      Consequently, \(\zeta^-=0\), and hence \(\zeta_1,\zeta_2\geq0\) in \((0,\infty)\).

      To obtain strict positivity, we write \eqref{eq:zeta}
      componentwise as
      \begin{align}
            -\Del\zeta_1 + V_{11}\zeta_1 + V_{12}\zeta_2 &= \frac{1}{2}f(1-f^2),\\
            -\Del\zeta_2 + V_{22}\zeta_2 + V_{12}\zeta_1 &= 0.
      \end{align}
      If \(\zeta_1(r_0)=0\) for some \(r_0>0\), then \(-\Del\zeta_1(r_0) \leq 0\) and \(V_{12}(r_0)\zeta_2(r_0) \leq 0\).
      Hence the first equation gives
      \begin{align}
            0 < \frac{1}{2}f(r_0)(1-f(r_0)^2) \leq 0,
      \end{align}
      where the strict inequality follows from \(0<f<1\) by Lemma~\ref{lem:shape}.
      This is a contradiction.      
      If \(\zeta_2(r_1)=0\) for some \(r_1>0\),
      then \(-\Del\zeta_2(r_1)\leq0\).
      Since \(V_{12}(r_1)<0\) and \(\zeta_1(r_1)>0\),
      the second equation gives
      \begin{align}
            0 = -\Del\zeta_2(r_1) + V_{12}(r_1)\zeta_1(r_1) < 0.
      \end{align}
      This is a contradiction.
      Thus \(\zeta_2>0\) in \((0,\infty)\).
      Since \(\partial_\lambda a_\lambda = (r/n)\zeta_2\), the inequalities \eqref{eq:monotone} follow.
\end{proof}

\section*{Acknowledgments}
\noindent
The author would like to thank Norihisa Ikoma for their valuable guidance and assistance in preparing and revising the manuscript.
This work was supported by JST SPRING, Japan Grant Number JPMJSP2123.

\noindent\textbf{AI assistance.}
The main results of this paper were initially obtained by Claude Opus 5.5.
The author subsequently verified all mathematical arguments and claims, simplified the proofs, and rewrote the entire manuscript.
ChatGPT was used to assist with English-language editing and stylistic improvements.
The author takes full responsibility for the mathematical correctness and all other aspects of this paper.

\bibliographystyle{unsrt}
\bibliography{ref}
\end{document}